\documentclass[12pt,english]{article}
\usepackage[latin9]{inputenc}
\usepackage{geometry}
\usepackage{babel}
\usepackage{amsmath}
\usepackage{amsthm}
\usepackage{amssymb}

\usepackage[unicode=true,pdfusetitle,
 bookmarks=true,bookmarksnumbered=false,bookmarksopen=false,
 breaklinks=false,pdfborder={0 0 1},backref=false,colorlinks=false]
 {hyperref}

\makeatletter
\theoremstyle{plain}
\newtheorem{thm}{\protect\theoremname}
  \theoremstyle{plain}
  \newtheorem{prop}[thm]{\protect\propositionname}
  \theoremstyle{plain}
  \newtheorem{lem}[thm]{\protect\lemmaname}
  \theoremstyle{plain}
  \newtheorem{cor}[thm]{\protect\corollaryname}
  \theoremstyle{plain}
  \newtheorem{claim}[thm]{\protect\claimname}
  \theoremstyle{plain}
  \newtheorem{case}[thm]{\protect\casename}

\date{}

\makeatother
  \providecommand{\casename}{Case}
  \providecommand{\claimname}{Claim}
  \providecommand{\corollaryname}{Corollary}
  \providecommand{\lemmaname}{Lemma}
  \providecommand{\propositionname}{Proposition}
\providecommand{\theoremname}{Theorem}

\begin{document}

\title{Uniqueness in Harper's vertex-isoperimetric theorem}

\author{Eero R\"{a}ty\thanks{Centre for Mathematical Sciences, Wilberforce Road, Cambridge CB3
0WB, UK, epjr2@cam.ac.uk}}
\maketitle
\begin{abstract}
For a set $A\subseteq Q_{n}=\left\{ 0,1\right\} ^{n}$ the $t$-neighbourhood
of $A$ is $N^{t}\left(A\right)=\left\{ x\,:\,d\left(x,A\right)\leq t\right\}$,
where $d$ denotes the usual graph distance on $Q_{n}$. Harper's
vertex-isoperimetric theorem states that among the subsets $A\subseteq Q_{n}$
of given size, the size of the $t$-neighbourhood is minimised when
$A$ is taken to be  an initial segment of the simplicial order. Aubrun
and Szarek asked the following question: if $A\subseteq Q_{n}$ is
a subset of given size for which the sizes of both $N^{t}\left(A\right)$
and $N^{t}\left(A^{c}\right)$ are minimal for all $t>0$, does it
follow that $A$ is isomorphic to an initial segment of the simplicial
order?

Our aim is to give a counterexample. Surprisingly it turns out that
there is no counterexample that is a Hamming ball, meaning a set that
lies between two consecutive exact Hamming balls, i.e.\ a set $A$
with $B\left(x,r\right)\subseteq A\subseteq B\left(x,r+1\right)$
for some $x\in Q_{n}$. We go further to classify all the sets $A\subseteq Q_{n}$
for which the sizes of both $N^{t}\left(A\right)$ and $N^{t}\left(A^{c}\right)$
are minimal for all $t>0$ among the subsets of $Q_{n}$ of given
size. We also prove that, perhaps surprisingly, if $A\subseteq Q_{n}$
for which the sizes of $N\left(A\right)$ and $N\left(A^{c}\right)$
are minimal among the subsets of $Q_{n}$ of given size, then the
sizes of both $N^{t}\left(A\right)$ and $N^{t}\left(A^{c}\right)$
are also minimal for all $t>0$ among the subsets of $Q_{n}$ of given
size. Hence the same classification also holds when we only require
$N\left(A\right)$ and $N\left(A^{c}\right)$ to have minimal size
among the subsets $A\subseteq Q_{n}$ of given size. 
\end{abstract}

\section{Introduction}

The $n$-dimensional hypercube $Q_{n}$ has vertex-set the power set
${\cal P}\left(\left\{ 1,\dots,n\right\} \right)$ with metric $d\left(x,y\right)=\left|x\Delta y\right|$.
For a subset $A$ of the hypercube $Q_{n}$, define the \textit{neighbourhood}
of $A$ to be the set $N\left(A\right)=\left\{ x\in Q_{n}:\,d\left(x,A\right)\leq1\right\} $
where $d\left(x,A\right)=\min_{y\in A}d\left(x,y\right)$. Also more
generally for each $t>0$ define $N^{t}\left(A\right)=\left\{ x\in Q_{n}:\,d(x,A)\leq t\right\} $,
and note that we have $N^{t}\left(A\right)=N\left(N^{t-1}\left(A\right)\right)$.

In order to state Harper's vertex-isoperimetric theorem we need a
few definitions. For any $n$ and $0\leq r\leq n$ define the \textit{lexicographic
order} on $\left\{ x:\,x\subseteq\left\{ 1,\dots,n\right\} ,\,|x|=r\right\} $
to be given by $x<_{lex}y$ if $\min\left(x\Delta y\right)\in x$,
and define the \textit{colexicographic order} to be given by $x<_{colex}y$
if $\text{max}\left(x\Delta y\right)\in y$. Define the \textit{simplicial
order} on $Q_{n}$ to be given by $x<_{sim}y$ if 

\[
|x|<|y|\text{ or }\left(|x|=|y|\text{ and }x<_{lex}y\right).
\]

\begin{thm}[Harper, \cite{key-13}]
Let $A$ be a
subset of $Q_{n}$ and let $B$ be the initial segment of the simplicial
order with $\left|A\right|=\left|B\right|$. Then $\left|N\left(A\right)\right|\geq\left|N\left(B\right)\right|$.
$\hfill\square$
\end{thm}
It turns out that the sets for which Harper's theorem holds with equality
are not in general unique. As a trivial example, any subset of $Q_{2}$
of size $2$ has minimal vertex boundary and not all such sets are
isomorphic. There are more interesting and less trivial examples as
well.

It is easy to verify that if $A$ is an initial segment of the simplicial
order, then so is $N\left(A\right)$. Hence Harper's theorem implies
that among the subsets $A\subseteq Q_{n}$ of given size, the size
of $N^{t}\left(A\right)$ is minimised when $A$ is chosen to be an
initial segment of the simplicial order. We say that $N^{t}\left(A\right)$
is \textit{minimal} if $\left|N^{t}\left(B\right)\right|\geq\left|N^{t}\left(A\right)\right|$
for all $B\subseteq Q_{n}$ with $\left|B\right|=\left|A\right|$.
Let $C$ be the initial segment of the simplicial order of size $\left|A\right|$.
It is useful to observe that Harper's theorem implies that $A$ is
minimal if and only if $\left|N^{t}\left(A\right)\right|=\left|N^{t}\left(C\right)\right|$. 

For a general introduction to the vertex-isoperimetric theorem, see
e.g.\ Bollob\'{a}s \cite[Chapter 16]{key-10}.

We say that two subsets $A$ and $B$ of $Q_{n}$ are \textit{isomorphic}
if there exists an isometry $\theta$ of $Q_{n}$ for which $\theta\left(A\right)=B$.
In this paper we consider the following question of Aubrun and Szarek
\cite[Exercise 5.66]{key-9}: If $A\subseteq Q_{n}$ for which $N^{t}\left(A\right)$
and $N^{t}\left(A^{c}\right)$ are minimal for all $t>0$, does it
follow that $A$ is isomorphic to an  initial segment of the simplicial
order? For convenience, we say that $A$ is \textit{extremal }if $N^{t}\left(A\right)$
and $N^{t}\left(A^{c}\right)$ are minimal for all $t>0$. Let $C$
be the initial segment of the simplicial order of size $\left|A\right|$.
Since $C^{c}$ is isomorphic to the initial segment of the simplicial
order of size $\left|C^{c}\right|$, it follows from Harper's inequality
that $A$ is extremal if and only if $\left|N^{t}\left(A\right)\right|=\left|N^{t}\left(C\right)\right|$
and $\left|N^{t}\left(A^{c}\right)\right|=\left|N^{t}\left(C^{c}\right)\right|$
for all $t$. In particular, initial segments of the simplicial order
are extremal. 

Define the \textit{exact Hamming ball} of radius $r$ centred at $x$
to be $B\left(x,r\right)=\left\{ y\in Q_{n}:\,d\left(x,y\right)\leq r\right\} $,
and define $A$ to be a \textit{Hamming ball} if there exist $x$
and $r$ such that $B\left(x,r\right)\subseteq A\subset B\left(x,r+1\right)$.
Note that $B\left(\emptyset,r\right)$ is the initial segment of the
simplicial order of size  $\sum_{i=0}^{r}{n \choose i}$, and every
initial segment of the simplicial order is a Hamming ball. Later in
the paper we sometimes consider exact Hamming balls of radius $r$
centred at $x$ on ${\cal P}\left(X_{i}\right)$ rather than on ${\cal P}\left(X\right)$.
In order to highlight this difference, we use $B_{i}\left(x,r\right)$
to denote the exact Hamming ball of radius $r$ centred at $x$ with
respect to ground set $X_{i}$. 

Note that requiring only $N^{t}\left(A\right)$ to be minimal for
all $t>0$ is not a strong enough condition to guarantee that $A$
should be isomorphic to an initial segment of the simplicial order.
Indeed, one could take for example $A=B\left(x,r\right)\setminus\left\{ x\right\} $
for $r\geq1$. Then $N^{t}\left(A\right)=B(x,r+t)$ for all $t>0$,
and hence $N^{t}\left(A\right)$ is always minimal, yet $A$ is not
isomorphic to an initial segment of the simplicial order. 

It turns out that the answer to the question is negative, and we present
a counterexample in Section 2. It turns out that all the extremal
sets are contained between two exact Hamming balls with same centre
and radius differing by 2 - that is, if $A$ is extremal then there
exist $x$ and $r$ with $B\left(x,r\right)\subseteq A\subset B\left(x,r+2\right)$.
Rather surprisingly, it turns out that the only Hamming balls which
are extremal \textit{are} the initial segments of the simplicial order. 

The second aim of this paper is to classify all the extremal sets
$A$ up to isomorphism. In order to state the result, we need some
notation. We write $X=\left\{ 1,\dots,n\right\} $, $X^{(r)}=\left\{ x\subseteq X\,:\,\left|x\right|=r\right\} $,
$X^{(\geq r)}=\left\{ x\subseteq X\,:\,\left|x\right|\geq r\right\} $,
$X_{i}=\left\{ 1,\dots,n\right\} \setminus\left\{ i\right\} $ and
$X_{i,j}=\left\{ 1,\dots,n\right\} \setminus\left\{ i,j\right\} $.
Throughout the paper, we denote the elements of $Q_{n}$ by lower
case letters, the subsets of $Q_{n}$ by upper case letters and the
set systems on $X^{\left(r\right)}$ by calligraphy letters. 

Define the maps $\pi_{i}\,:\,X^{\left(r+1\right)}\rightarrow X_{i}^{\left(r\right)}\cup X_{i}^{\left(r+1\right)}$\textit{
}by $\pi_{i}\left(x\right)=x\setminus\left\{ i\right\} $ for all
$x\in X^{\left(r+1\right)}$. For a set system ${\cal B}\subseteq X^{\left(r+1\right)}$
define $\pi_{i}\left({\cal B}\right)=\left\{ \pi_{i}\left(x\right)\,:\,x\in{\cal B}\right\} $.
Note that $\pi_{i}$ is a bijection from $X^{\left(r+1\right)}$ to
$X_{i}^{\left(r\right)}\cup X_{i}^{\left(r+1\right)}$, and hence
$\left|\pi_{i}\left({\cal B}\right)\right|=\left|{\cal B}\right|$
for all $i$. 

It is known that the exact Hamming balls are the only sets of their
respective sizes for which the inequality in Harper's theorem holds
with equality. That is, if $A\subseteq Q_{n}$ is a set of size $\left|X^{(\leq r)}\right|$
for which the size of $N\left(A\right)$ is minimal among the subsets
of $Q_{n}$ of the same size, then $A=B\left(x,r\right)$ for some
$x\in Q_{n}$. Thus if $A\subseteq Q_{n}$ is a set of size $\left|X^{(\leq r)}\right|$
for some $r$ and $A$ is extremal, it certainly follows that $A$
is isomorphic to the initial segment of the simplicial order. 

Note that if $A$ is extremal then so is $A^{c}$, as the definition
of extremality is symmetric under taking complements. Set $G_{r}=X^{(\leq r)}\cup\left\{ B\in X^{(r+1)}\,:\,1\in B\right\} $.
It is easy to check that $\left|G_{r}\right|+\left|G_{n-r-2}\right|=2^{n}$.
Thus if $A\subseteq Q_{n}$ satisfies $\left|A\right|\neq\left|X^{(\leq r)}\right|$
for all $r$, then at least one of $\left|X^{(\leq r)}\right|<\left|A\right|\leq\left|G_{r}\right|$
or $\left|X^{(\leq r)}\right|<\left|A^{c}\right|\leq\left|G_{r}\right|$
is satisfied for some $r$. Hence it is sufficient to classify only
those extremal sets $A\subseteq Q_{n}$ for which there exists $r$
such that $\left|X^{\left(\leq r\right)}\right|<\left|A\right|\leq\left|G_{r}\right|$. 

For convenience, we write $f_{r}=f_{n,r}=\left|X^{(\le r)}\right|=\sum_{j=0}^{r}{n \choose j}$
and $g_{r}=g_{n,r}=\left|G_{r}\right|=\sum_{j=0}^{r}{n \choose j}+{n-1 \choose r}$.
In both cases the dependence on $n$ will not be highlighted if the
value of $n$ is clear from the context.

Let $s$ be an integer of the form $s=f_{r}+k$ for some $0\leq k\leq{n-1 \choose r}$.
Note that for fixed $n$, the value of $s$ uniquely determines the
values of $r$ and $k$. Furthermore, observe that $f_{n}+{n-1 \choose r}=g_{r}$.
Hence the set of integers that can be written on that form is exactly
the set of those integers $s'$ for which there exists $r$ such that
$f_{r}\leq s'\leq g_{r}$. 

Given an integer $s$ of the form $s=f_{r}+k$, let ${\cal A}$ be
the initial segment of the lexicographic order on $X^{\left(r+1\right)}$
of size  $k$. For each $i$ set 
\[
A_{i}=X^{\left(\leq r\right)}\cup\left(\left\{ i\right\} +\pi_{i}\left({\cal A}\right)\right).
\]
Note that $\pi_{i}\left({\cal A}\right)\subseteq X_{i}^{\left(r\right)}\cup X_{i}^{\left(r+1\right)}$,
so $\left\{ i\right\} +\pi_{i}\left({\cal A}\right)$ is a well-defined
subset of $X^{\left(r+1\right)}\cup X^{\left(r+2\right)}$ of the
same size as ${\cal A}$. Also note that the sets $X^{\left(\leq r\right)}$
and $\left\{ i\right\} +\pi_{i}\left({\cal A}\right)$ are disjoint.
Hence each $A_{i}$ has size $s$. Furthermore, since $k\leq{n-1 \choose r}$,
it follows that ${\cal A}\subseteq\left\{ 1\right\} +X_{1}^{\left(r\right)}$.
Hence $\left\{ 1\right\} +\pi_{1}\left({\cal A}\right)={\cal A}$,
so $A_{1}$ is just the initial segment of the simplicial order. Note
that some of the $A_{i}$'s might be isomorphic to each other. 

Now we are ready to state the classification of extremal sets.
\begin{thm}[Classification of extremal sets]
\label{thm:Theorem 2}
Let $A\subseteq Q_{n}$ be a subset of size $s$, where $s=f_{r}+k$
for some $r$ and $k\leq{n-1 \choose r}$. Let $A_{1},\dots,A_{n}$
be the sets defined as above for these choices of $r$ and $k$. Then
$A$ is extremal if and only if $A$ is isomorphic to some $A_{i}$. 
\end{thm}
It is natural to ask what happens if we weaken the notion of extremality.
A natural way to do this is to seek for subsets $A\subseteq Q_{n}$
for which both $N\left(A\right)$ and $N\left(A^{c}\right)$ have
minimal size among the subsets of $Q_{n}$ of given size. We say that
$A\subseteq Q_{n}$ is \textit{weakly extremal} if for all $B\subseteq Q_{n}$
with $\left|B\right|=\left|A\right|$ we have $\left|N\left(B\right)\right|\geq\left|N\left(A\right)\right|$
and $\left|N\left(B^{c}\right)\right|\geq\left|N\left(A^{c}\right)\right|$.
Let $C$ be the initial segment of the simplicial order of size  $\left|A\right|$.
Again by Harper's theorem, weak extremality of $A$ is equivalent
to having $\left|N\left(A\right)\right|=\left|N\left(C\right)\right|$
and $\left|N\left(A^{c}\right)\right|=\left|N\left(C^{c}\right)\right|$.
Rather surprisingly, we prove in Section 4 that the notions of weak
extremality and extremality coincide. 
\begin{thm}
\label{thm:Theorem 3}Let $A\subseteq Q_{n}$ be a subset for which
every $B\subseteq Q_{n}$ with $\left|B\right|=\left|A\right|$ satisfies
$\left|N\left(B\right)\right|\geq\left|N\left(A\right)\right|$ and
$\left|N\left(B^{c}\right)\right|\geq\left|N\left(A^{c}\right)\right|$.
Then for all $t>0$ and $B\subseteq Q_{n}$ with $\left|B\right|=\left|A\right|$
we have $\left|N^{t}\left(B\right)\right|\geq\left|N^{t}\left(A\right)\right|$
and $\left|N^{t}\left(B^{c}\right)\right|\geq\left|N^{t}\left(A^{c}\right)\right|$.
\end{thm}
Thus it immediately follows that Theorem \ref{thm:Theorem 2} holds
when extremality is replaced by weak extremality. The main ingredient
in the proof of Theorem \ref{thm:Theorem 3}  is the following theorem,
which we also prove in Section 4. 
\begin{thm}
\label{thm:Theorem 4}Let $A\subseteq Q_{n}$ be a subset
for which every $B\subseteq Q_{n}$ with $\left|B\right|=\left|A\right|$
satisfies $\left|N\left(B\right)\right|\geq\left|N\left(A\right)\right|$.
Then for all $t>0$ and $B\subseteq Q_{n}$ with $\left|B\right|=\left|A\right|$
we also have $\left|N^{t}\left(B\right)\right|\geq\left|N^{t}\left(A\right)\right|$.
\end{thm}
The plan of the paper is as follows. In Section 2 we construct an
extremal set which is not isomorphic to an initial segment of the
simplicial order. In Section 3 we prove Theorem \ref{thm:Theorem 2}.
In Section 4 we consider weakly extremal sets and prove Theorem \ref{thm:Theorem 3}. 

Recall that exact Hamming balls are the unique sets of their respective
sizes for which equality holds in Harper's inequality. In Section
5 we prove another near-uniqueness result: we show that there exists
only one set $B_{r}$ of size $g_{r}$, apart from the initial segment,
for which equality holds in Harper's inequality. In fact the set $B_{r}$
is also an extremal set, and we describe it already in Section 2.

\section{Construction of an example}

In this section we find for every $r$ an extremal set $B_{r}\subseteq Q_{n}$
satisfying $\left|B_{r}\right|=g_{r}$ for which $B_{r}$ is not isomorphic
to the initial segment of the simplicial order. Let $C_{r}$ be the
initial segment of the simplicial order of size $g_{r}$ on $Q_{n}$.
Then $C_{r}$ can be written as $C_{r}=B\left(\emptyset,r\right)\cup B\left(\left\{ 1\right\} ,r\right)$.
Define $B_{r}=B\left(\emptyset,r\right)\cup B\left(\left\{ 1,2\right\} ,r\right)$.
In particular, $B_{r}$ is the union of two exact Hamming balls with
same radius and centres within distance 2 apart from each other. 

Note that if $B\left(s,r\right)\subseteq B_{r}\subseteq B\left(s,r+1\right)$
holds for some $s\in Q_{n}$, then the first inclusion implies that
$s=\emptyset$ or $s=\left\{ 1,2\right\} $. However, the second inclusion
is violated in both cases. Thus $B_{r}$ is not a Hamming ball, and
hence it is not isomorphic to the initial segment of the simplicial
order. 

It is easy to verify that $N^{t}\left(C_{r}\right)=C_{r+t}$ and $N^{t}\left(B_{r}\right)=B_{r+t}$
for all $t>0$, so in order to prove that $N^{t}\left(B_{r}\right)$
are minimal for all $t>0$ it suffices to prove that $\left|B_{m}\right|=\left|C_{m}\right|$
for all $m$. Observe that $B_{r}$ can be written as $B_{r}=X^{\left(\leq r\right)}\cup\left(\left\{ 1,2\right\} +\left(X_{1,2}^{\left(r-1\right)}\cup X_{1,2}^{\left(r\right)}\right)\right)$.
Hence 
\[
\left|B_{r}\right|=f_{r}+{n-2 \choose r-1}+{n-2 \choose r}=f_{r}+{n-1 \choose r}=g_{r}
\]
and thus $\left|B_{r}\right|=\left|C_{r}\right|$ for all $r$, as
required. 

Also note that
\[
C_{r}^{c}=B\left(\left\{ 1,\dots,n\right\} ,n-r-1\right)\cap B\left(\left\{ 2,\dots,n\right\} ,n-r-1\right)
\]
and hence it is easy to check that 
\[
C_{r}^{c}=B\left(\left\{ 1,\dots,n\right\} ,n-r-2\right)\cup B\left(\left\{ 2,\dots,n\right\} ,n-r-2\right).
\]
Thus $C_{r}^{c}$ is isomorphic to $C_{n-r-2}$ under the isometry
$\theta\,:\,Q_{n}\rightarrow Q_{n}$ given by $\theta\left(a\right)=a^{c}$
for all $a\in Q_{n}$. In particular, it follows that $g_{r}+g_{n-r-2}=2^{n}$. 

Similarly we have 
\[
B_{r}^{c}=B\left(\left\{ 1,\dots,n\right\} ,n-r-1\right)\cap B\left(\left\{ 3,\dots,n\right\} ,n-r-1\right)
\]
and our aim is to show that this implies that 
\begin{equation}
B_{r}^{c}=B\left(\left\{ 1,3,\dots,n\right\} ,n-r-2\right)\cup B\left(\left\{ 2,\dots,n\right\} ,n-r-2\right).\label{eq:1}
\end{equation}
Indeed, note that for all $x\in B\left(\left\{ 1,3,\dots,n\right\} ,n-r-2\right)$
we have 
\[
d\left(x,\left\{ 1,\dots,n\right\} \right)\leq d\left(x,\left\{ 1,3\dots,n\right\} \right)+d\left(\left\{ 1,3,\dots,n\right\} ,\left\{ 1,\dots,n\right\} \right)\leq\left(n-r-2\right)+1=n-r-1
\]
by the triangle inequality. One can similarly deduce that $d\left(x,\left\{ 3,\dots,n\right\} \right)\leq n-r-1$
holds by the triangle inequality. Hence 
\[
B\left(\left\{ 1,3,\dots,n\right\} ,n-r-2\right)\subseteq B\left(\left\{ 1,\dots,n\right\} ,n-r-1\right)\cap B\left(\left\{ 3,\dots,n\right\} ,n-r-1\right),
\]
and by using the same argument it follows that 
\[
B\left(\left\{ 2,\dots,n\right\} ,n-r-2\right)\subseteq B\left(\left\{ 1,\dots,n\right\} ,n-r-1\right)\cap B\left(\left\{ 3,\dots,n\right\} ,n-r-1\right).
\]
These two observations together imply the $\supseteq$-part of (\ref{eq:1}). 

Note that $B\left(\left\{ 1,3,\dots,n\right\} ,n-r-2\right)\cup B\left(\left\{ 2,\dots,n\right\} ,n-r-2\right)$
is isomorphic to $B_{n-r-2}$ under the isometry $\phi$ given by
$\phi\left(a\right)=a\Delta\left\{ 2,\dots,n\right\} $ for all $a\in Q_{n}$.
Hence 
\[
\left|B\left(\left\{ 1,3\dots,n\right\} ,n-r-2\right)\cup B\left(\left\{ 2,\dots,n\right\} ,n-r-2\right)\right|=\left|B_{n-r-2}\right|=g_{n-r-2}=2^{n}-\left|B_{r}\right|=\left|B_{r}^{c}\right|.
\]
This together with the fact that the inclusion holds in the $\supseteq$-direction
completes the proof of (\ref{eq:1}). In particular, $B_{r}^{c}$
is isomorphic to $B_{n-r-2}$. 

Hence it follows that $N^{t}\left(B_{r}^{c}\right)$ is isomorphic
to $B_{n-r+t-2}$ and $N^{t}\left(C_{r}^{c}\right)$ is isomorphic
to $C_{n-r+t-2}$ for all $t>0$. Since $\left|B_{n-r+t-2}\right|=\left|C_{n-r+t-2}\right|$,
it follows that $\left|N^{t}\left(B_{r}^{c}\right)\right|=\left|N^{t}\left(C_{r}^{c}\right)\right|$
for all $t>0$ and hence $N^{t}\left(B_{r}^{c}\right)$ are minimal
for all $t>0$. Thus $B_{r}$ is an extremal set. 

\section{Classifying all extremal sets}

Recall that $f_{r}=\sum_{i=0}^{r}{n \choose i}$ is the size of the
exact Hamming ball of radius $r$, and $g_{r}=\sum_{i=0}^{r}{n \choose i}+{n-1 \choose r}$
is the size of the initial segment $X^{(\leq r)}\cup\left(\left\{ 1\right\} +X_{1}^{(r)}\right)$,
where $X_{i}=\left\{ 1,\dots,n\right\} \setminus\left\{ i\right\} $.
It is convenient to exclude sets of size $f_{r}$ from the classification,
and this is possible due to the following much stronger result. 
\begin{prop}
\label{prop:Prop5}Let $A\subseteq Q_{n}$ be a set satisfying $\left|A\right|=f_{r}$,
and such that for any $B\subseteq Q_{n}$ with $\left|B\right|=f_{r}$
we have $\left|N\left(B\right)\right|\geq\left|N\left(A\right)\right|$.
Then $A=B\left(x,r\right)$ for some $x\in Q_{n}$. $\hfill\square$
\end{prop}
Since this is a well-known fact, the proof is omitted. This could
be deduced by induction on $n$ and applying Lemma 6 of Katona from
\cite{key-14}. A similar technique will be used in the proof of Claim
1 in Theorem \ref{thm:Theorem 21} in Section 5.

Since the classification of extremal sets $A\subseteq Q_{n}$ satisfying
 $\left|A\right|=f_{r}$ for some $r$ is covered by Proposition \ref{prop:Prop5},
it is enough to consider only those sets $A\subseteq Q_{n}$ satisfying
$f_{r}<\left|A\right|<f_{r+1}$ for some $r$. Furthermore, since
$g_{r}+g_{n-2-r}=2^{n}$ and $f_{r}+f_{n-1-r}=2^{n}$, by considering
$A^{c}$ if necessary, it is enough to classify only those extremal
sets $A\subseteq Q_{n}$ satisfying $f_{r}<\left|A\right|\leq g_{r}$
for some $r$. Hence from now on we will assume that $A\subseteq Q_{n}$
is an extremal set satisfying $f_{r}<\left|A\right|\leq g_{r}$ for
some $r$.
\begin{lem}
\label{lem:Lemma6}Let $A\subseteq Q_{n}$ be an extremal set which
satisfies $f_{r}<\left|A\right|\leq g_{r}$ for some $r$. Then there
exist $x,\,y,\,z\in Q_{n}$ with $y\neq z$ satisfying $d\left(x,y\right)\leq1$,
$d\left(x,z\right)\leq1$ and $B\left(x,r\right)\subseteq A\subseteq B\left(y,r+1\right)\cap B\left(z,r+1\right)$. 
\end{lem}
Note that since $d\left(x,y\right)\leq1$, the condition $A\subseteq B\left(y,r+1\right)$
implies that $A\subseteq B\left(x,r+2\right)$. Hence it follows that
there exists $x\in Q_{n}$ with $B\left(x,r\right)\subseteq A\subseteq B\left(x,r+2\right)$.
This implies that the interesting behaviour in the set $A$ occurs
only on two layers of the cube, namely on those which are at distance
$r+1$ and $r+2$ apart from $x$. 
\begin{proof}
Let $A\subseteq Q_{n}$ with $f_{r}<\left|A\right|\leq g_{r}$. Let
$C_{r}$ be the initial segment of size $g_{r}$, and recall from
Section 2 that $N^{t}\left(C_{r}\right)=C_{r+t}$. Since $\left|A\right|\leq g_{r}$
and $A$ is extremal, it follows that $\left|N^{n-r-2}\left(A\right)\right|\leq\left|N^{n-r-2}\left(C_{r}\right)\right|=g_{n-2}=2^{n}-2$.
Hence there exist distinct elements $u,\,v\in N^{n-r-2}\left(A\right)^{c}$.
Let $y=u^{c}$ and $z=v^{c}$. By the choice of $u$ and $v$ it follows
that $B\left(u,n-r-2\right)\cup B\left(v,n-r-2\right)\subseteq A^{c}$.
Taking complements gives that $A\subseteq B\left(y,r+1\right)\cap B\left(z,r+1\right)$. 

Since $\left|A\right|>f_{r}$, it follows that $\left|A^{c}\right|<f_{n-r-1}$.
Combining this with the extremality of $A$, it follows that $\left|N^{r}\left(A^{c}\right)\right|\leq\left|N^{r}\left(B\left(\emptyset,n-r-1\right)\right)\right|=f_{n-1}=2^{n}-1$.
Hence there exists $x\in\left(N^{r}\left(A^{c}\right)\right)^{c}$,
which implies that $B\left(x,r\right)\cap A^{c}=\emptyset$. Hence
it follows that $B\left(x,r\right)\subseteq A$, and combining this
with the previous observations implies that $B\left(x,r\right)\subseteq A\subseteq B\left(y,r+1\right)\cap B\left(z,r+1\right)$.
Since $B\left(x,r\right)\subseteq B\left(y,r+1\right)$, we must have
$d\left(x,y\right)\leq1$ and similarly it follows that $d\left(x,z\right)\leq1$,
as required. 
\end{proof}
The proof of Lemma \ref{lem:Lemma6} gives some insight on why it
is convenient to assume that the size of $A$ satisfies $f_{r}<\left|A\right|\leq g_{r}$
rather than just $f_{r}<\left|A\right|<f_{r+1}$. Indeed, the condition
$f_{r}<\left|A\right|<f_{r+1}$ would not be strong enough to guarantee
the existence of both $y$ and $z$.

Given this result, we can split the rest of the classification into
two parts: considering those $A$ which are Hamming balls, i.e.\ for
which there exist $x\in Q_{n}$ and $r$ satisfying $B\left(x,r\right)\subseteq A\subseteq B\left(x,r+1\right)$,
and considering those $A$ for which no such $x$ and $r$ exist.
It turns out that all the examples apart from the initial segment
appear in the second case. This is proved in Proposition \ref{prop:Prop 11},
but before that we need a few preliminary results. Many of these preliminary
results are used later as well. 

For $A\subseteq Q_{n}$, define the $i$\textit{-sections }$A_{i,\,+}$
and $A_{i,\,-}$ of $A$ by setting $A_{i,\,+}=\left\{ a\setminus\left\{ i\right\} \,:\,a\in A,\,i\in a\right\} $
and $A_{i,\,-}=\left\{ a\,:\,a\in A,\,i\not\in a\right\} $. Note
that these are subsets of ${\cal P}\left(X_{i}\right)$, which is
naturally isomorphic to $Q_{n-1}$. If the direction $i$ is clear
from the context, these will be denoted as $A_{+}$ and $A_{-}$.

We often want to relate the neighbourhood of $A$ to the neighbourhoods
of the $i$-sections $A_{i,\,+}$ and $A_{i,\,-}$. However, for $A_{i,\,+}$
and $A_{i,\,-}$ the neighbourhood is always taken inside ${\cal P}\left(X_{i}\right)$,
i.e.\ with respect to the ground set $X_{i}$ rather than $X$. Since
for $i$-sections $A_{i,\,+}$ and $A_{i,\,-}$ we always consider
the neighbourhood with respect to $X_{i}$, and otherwise we always
consider neighbourhood with respect to $X$, we will use the same
notation $N\left(A\right)$ and $N\left(A_{i,\,\pm}\right)$ in both
cases to avoid excessive use of subscripts in our notation. In the
second case the neighbourhood $N\left(A_{i,\,\pm}\right)$ that is
considered is a neighbourhood of an $i$-section, and hence it should
always be understood as the neighbourhood with respect to the ground
set $X_{i}$. 
\begin{lem}
\label{lem:Lemma 7}For all $r\geq0$ and for any distinct elements
$x,\,y\in Q_{n}$ we have $\left|B\left(x,r\right)\cup B\left(y,r\right)\right|\geq g_{r}$,
with equality if and only if $d\left(x,y\right)\leq2$. 
\end{lem}
\begin{proof}
By the symmetry of $Q_{n}$ we may assume that $x=\emptyset$. Set
$A=B\left(\emptyset,r\right)\cup B\left(y,r\right)$. Recall that
$B_{i}\left(x,r\right)$ denotes the exact Hamming ball of radius
$r$ centred at $x$ with respect to the ground set $X_{i}$. For
any $i\in y$ we have $A_{i,\,-}=B_{i}\left(\emptyset,r\right)\cup B_{i}\left(y\setminus\left\{ i\right\} ,r-1\right)$
and $A_{i,\,+}=B_{i}\left(\emptyset,r-1\right)\cup B_{i}\left(y\setminus\left\{ i\right\} ,r\right)$.
In particular, we have $B_{i}\left(\emptyset,r\right)\subseteq A_{i,\,-}$
and $B_{i}\left(y\setminus\left\{ i\right\} ,r\right)\subseteq A_{i,\,+}$,
so $\left|A\right|=\left|A_{i,\,+}\right|+\left|A_{i,\,-}\right|\geq2f_{n-1,\,r}=g_{r}$,
which proves the first part. 

Note that the equality holds if and only if $A_{i,\,-}=B_{i}\left(\emptyset,r\right)$
and $A_{i,\,+}=B_{i}\left(y\setminus\left\{ i\right\} ,r\right)$.
Hence we must have $B_{i}\left(y\setminus\left\{ i\right\} ,r-1\right)\subseteq B_{i}\left(\emptyset,r\right)$
and $B_{i}\left(\emptyset,r-1\right)\subseteq B_{i}\left(y\setminus\left\{ i\right\} ,r\right)$,
and hence it follows that $d\left(y\setminus\left\{ i\right\} ,\emptyset\right)=d\left(y\setminus\left\{ i\right\} ,x\right)\leq1$.
Since $i\in y$, this holds if and only if $d\left(x,y\right)=d\left(y,\emptyset\right)\le2$,
which completes the proof of Lemma \ref{lem:Lemma 7}. 
\end{proof}
\begin{lem}
\label{lem:Lemma 8}Let $G$ be the set of isometries $\phi$ of $Q_{n}$
satisfying $\phi\left(\emptyset\right)=\emptyset$. For each $\sigma\in S_{n}$,
define $\phi_{\sigma}$ by setting $\phi_{\sigma}\left(a\right)=\left\{ \sigma\left(i\right)\,:\,i\in a\right\} $
for $a\in Q_{n}$. Then $G=\left\{ \phi_{\sigma}\,:\,\sigma\in S_{n}\right\} $. 
\end{lem}
\begin{proof}
It is clear that each $\phi_{\sigma}$ is an isometry and $\phi_{\sigma}\left(\emptyset\right)=\emptyset$.
Let $g\in G$. Since $g\left(\emptyset\right)=\emptyset$, it follows
that for all $i\in X$ we have $d\left(g\left(\left\{ i\right\} \right),\emptyset\right)=1$.
Hence $g\left(\left\{ i\right\} \right)$ is a set with exactly one
element. Let this element be $g_{i}$, and define $\sigma$ by setting
$\sigma\left(i\right)=g_{i}$. Note that since $g$ is an injection,
it follows that $\sigma$ is an injection and hence $\sigma\in S_{n}$.

Our aim is to prove that $g\left(a\right)=\phi_{\sigma}\left(a\right)$
for all $a\in Q_{n}$. We prove this by induction on the number of
elements in $a$. Note that the claim is true for any $a\in Q_{n}$
with $\left|a\right|\in\left\{ 0,1\right\} $ by the construction
of $\sigma$ and the fact that $g\left(\emptyset\right)=\emptyset$. 

Suppose that the claim is true for all $b\in Q_{n}$ with $\left|b\right|\leq m-1$,
where $m\geq2$, and let $a\in Q_{n}$ with $\left|a\right|=m$. Let
$i\in a$. Then $g\left(a\setminus\left\{ i\right\} \right)=\phi_{\sigma}\left(a\setminus\left\{ i\right\} \right)$
by induction. Note that $\left|\phi_{\sigma}\left(a\setminus\left\{ i\right\} \right)\right|=\left|a\right|-1=m-1$
and $\left|g\left(a\right)\right|=d\left(g\left(a\right),\emptyset\right)=d\left(a,\emptyset\right)=m$.
We also have $d\left(g\left(a\right),g\left(a\setminus\left\{ i\right\} \right)\right)=d\left(a,a\setminus\left\{ i\right\} \right)=1$.
In particular, $g\left(a\right)$ is a set with $m$ elements, $g\left(a\setminus\left\{ i\right\} \right)$
is a set with $m-1$ elements and the symmetric difference of $g\left(a\right)$
and $g\left(a\setminus\left\{ i\right\} \right)$ contains one element.
Hence $g\left(a\right)=g\left(a\setminus\left\{ i\right\} \right)\cup\left\{ j\right\} =\phi_{\sigma}\left(a\setminus\left\{ i\right\} \right)\cup\left\{ j\right\} $
for some $j\not\in g\left(a\setminus\left\{ i\right\} \right)$. 

Note that $d\left(g\left(a\right),g\left(\left\{ i\right\} \right)\right)=d\left(a,\left\{ i\right\} \right)=m-1$.
Combining this with $\left|g\left(a\right)\right|=m$ and $\left|g\left(\left\{ i\right\} \right)\right|=\left|\left\{ \sigma\left(i\right)\right\} \right|=1$,
it follows that $\sigma\left(i\right)\in g\left(a\right)$. But $\sigma\left(i\right)\not\in\phi_{\sigma}\left(a\setminus\left\{ i\right\} \right)$,
so we must have $j=\sigma\left(i\right)$. Hence $g\left(a\right)=\phi_{\sigma}\left(a\setminus\left\{ i\right\} \right)\cup\left\{ \sigma\left(i\right)\right\} =\phi_{\sigma}\left(a\right)$,
as required. Hence the result follows by induction. 
\end{proof}
Note that for all $\sigma\in S_{n}$, $\phi_{\sigma}$ induces a bijection
on $X^{\left(r\right)}$. We say that ${\cal A}\subseteq X^{\left(r\right)}$
and ${\cal B}\subseteq X^{\left(r\right)}$ are \textit{isomorphic
as subsets of }$X^{\left(r\right)}$ if there exists $\phi_{\sigma}$
such that $\phi_{\sigma}\left({\cal A}\right)={\cal B}$. In order
to avoid confusion between this notion and the notion of being isomorphic
on $Q_{n}$, we will explicitly write down ``isomorphic as subsets
of $X^{\left(r\right)}$'' rather than ``isomorphic''. However,
the aim of the next Lemma is to relate these two ways of being isomorphic
together. 
\begin{lem}
\label{lem:Lemma 9}Let $A$ and $C$ be subsets of $Q_{n}$ of the
form $A=B\left(\emptyset,r\right)\cup{\cal A}$ and $C=B\left(\emptyset,r\right)\cup{\cal C}$
for some ${\cal A},\,{\cal C}\subseteq X^{\left(r+1\right)}$. If
$A$ and $C$ are isomorphic, then ${\cal A}$ and ${\cal C}$ are
isomorphic as subsets of $X^{\left(r\right)}$. 
\end{lem}
\begin{proof}
Note that it suffices to show that if $A$ and $C$ are isomorphic,
then there exists $\phi_{\sigma}\in G$ with $\phi_{\sigma}\left(A\right)=C$.
Let $g$ be an isometry mapping $A$ to $C$. Hence $g$ also maps
$A^{c}$ to $C^{c}$. If $g\left(\emptyset\right)=\emptyset$ we are
done by Lemma \ref{lem:Lemma 8}, so we may assume that $g\left(\emptyset\right)\neq\emptyset$. 

Since $g$ is an isometry, it follows that $B\left(g\left(\emptyset\right),r\right)=g\left(B\left(\emptyset,r\right)\right)$,
and hence $B\left(g\left(\emptyset\right),r\right)\subseteq C$. Thus
$B\left(\emptyset,r\right)\cup B\left(g\left(\emptyset\right),r\right)\subseteq C$,
and by Lemma \ref{lem:Lemma 7} we have $\left|C\right|\geq g_{r}$.
Furthermore, the equality holds if and only if $C=B\left(\emptyset,r\right)\cup B\left(g\left(\emptyset\right),r\right)$
and $d\left(g\left(\emptyset\right),\emptyset\right)\leq2$. 

Since $B\left(\left\{ 1,\dots,n\right\} ,n-r-2\right)\subseteq A^{c}$,
it follows similarly that $B\left(g\left(\left\{ 1,\dots,n\right\} \right),n-r-2\right)\subseteq C^{c}$.
Since $g$ is an isometry with $g\left(\emptyset\right)\neq\emptyset$,
it follows that $g\left(\left\{ 1,\dots,n\right\} \right)\neq\left\{ 1,\dots,n\right\} $
as $\emptyset$ and $\left\{ 1,\dots,n\right\} $ are antipodal points.
Hence 
\[
B\left(\left\{ 1,\dots,n\right\} ,n-r-2\right)\cup B\left(g\left(\left\{ 1,\dots,n\right\} \right),n-r-2\right)\subseteq C^{c}
\]
and thus Lemma \ref{lem:Lemma 7} implies that $\left|C^{c}\right|\geq g_{n-r-2}$. 

Since $g_{r}+g_{n-r-2}=2^{n}$, the equality must hold in both cases.
Thus $C=B\left(\emptyset,r\right)\cup B\left(g\left(\emptyset\right),r\right)$,
and since $C\subseteq B\left(\emptyset,r+1\right)$ we must have $d\left(g\left(\emptyset\right),\emptyset\right)=1$.
Hence $A=B\left(g^{-1}\left(\emptyset\right),r\right)\cup B\left(\emptyset,r\right)$,
and we also have $d\left(g^{-1}\left(\emptyset\right),\emptyset\right)=1$
since $g$ is an isometry. 

Let $i,\,j\in X$ be chosen such that $g\left(\emptyset\right)=\left\{ i\right\} $
and $g^{-1}\left(\emptyset\right)=j$. Since $A=B\left(\left\{ j\right\} ,r\right)\cup B\left(\emptyset,r\right)$
and $C=B\left(\left\{ i\right\} ,r\right)\cup B\left(\emptyset,r\right)$,
by choosing $\sigma=\left(ij\right)\in S_{n}$ we have that $\phi_{\sigma}\left(A\right)=C$,
as required. 
\end{proof}
\begin{lem}
\label{lem:Lemma 10}Let $t$ be a positive integer, $A\subseteq Q_{n}$
and let $j$ be a direction for which $\left|A_{j,\,-}\right|\geq\left|A_{j,\,+}\right|$.
Let $C_{j,\,-}$ and $C_{j,\,+}$ be the initial segments of the simplicial
order on ${\cal P}\left(X_{j}\right)$ with $\left|C_{j,\,-}\right|=\left|A_{j,\,-}\right|$
and $\left|C_{j,\,+}\right|=\left|A_{j,\,+}\right|$. If $N^{t}\left(A\right)$
is minimal and $C_{j,\,-}\subseteq N\left(C_{j,\,+}\right)$, then
$N^{t}\left(A_{j,\,-}\right)$ and $N^{t}\left(A_{j,\,+}\right)$
are minimal, and we have $N^{t-1}\left(A_{j,\,+}\right)\subseteq N^{t}\left(A_{j,\,-}\right)$
and $N^{t-1}\left(A_{j,\,-}\right)\subseteq N^{t}\left(A_{j,\,+}\right)$.
\end{lem}
The idea of comparing the neighbourhoods of $A_{+}$ and $A_{-}$
with the neighbourhoods of $C_{+}$ and $C_{-}$ is used similarly
in Bollob\'{a}s and Leader \cite{key-11}. 
\begin{proof}
For simplicity denote the $j$-sections by $A_{+}$, $A_{-}$,
$C_{+}$ and $C_{-}$. Define $C\subseteq Q_{n}$ by setting $C=C_{-}\cup\left(C_{+}+\left\{ j\right\} \right)$.
It is easy to observe that $\left(N^{t}\left(A\right)\right)_{+}=N^{t-1}\left(A_{-}\right)\cup N^{t}\left(A_{+}\right)$,
and we similarly have $\left(N^{t}\left(A\right)\right)_{-}=N^{t-1}\left(A_{+}\right)\cup N^{t}\left(A_{-}\right)$.
Hence it follows that
\begin{equation}
\left|N^{t}\left(A\right)\right|=\left|N^{t-1}\left(A_{-}\right)\cup N^{t}\left(A_{+}\right)\right|+\left|N^{t-1}\left(A_{+}\right)\cup N^{t}\left(A_{-}\right)\right|.\label{eq:L2}
\end{equation}
Since $C_{+}$ and $C_{-}$ are initial segments of the simplicial
order on ${\cal P}\left(X_{j}\right)$ with $\left|C_{-}\right|\geq\left|C_{+}\right|$,
it follows that $C_{+}\subseteq C_{-}$ as initial segments are nested.
Hence we also have $N^{t-1}\left(C_{+}\right)\subseteq N^{t-1}\left(C_{-}\right)\subseteq N^{t}\left(C_{-}\right)$,
as taking neighbourhoods preserves inclusions, and a set is
always a subset of its neighbourhood. Similarly the assumption $C_{-}\subseteq N\left(C_{+}\right)$
implies that we have $N^{t-1}\left(C_{-}\right)\subseteq N^{t}\left(C_{+}\right)$.
Thus applying (\ref{eq:L2}) for $C$ implies that 
\begin{equation}
\left|N^{t}\left(C\right)\right|=\left|N^{t}\left(C_{+}\right)\right|+\left|N^{t}\left(C_{-}\right)\right|.\label{eq:L3}
\end{equation}

Since $C_{+}$ and $C_{-}$ are initial segments, Harper's inequality
implies that we have $\left|N^{t}\left(A_{+}\right)\right|\geq\left|N^{t}\left(C_{+}\right)\right|$
and $\left|N^{t}\left(A_{-}\right)\right|\geq\left|N^{t}\left(C_{-}\right)\right|$.
In particular, it follows that 
\begin{equation}
\left|N^{t-1}\left(A_{-}\right)\cup N^{t}\left(A_{+}\right)\right|\geq\left|N^{t}\left(A_{+}\right)\right|\geq\left|N^{t}\left(C_{+}\right)\right|\label{eq:L4}
\end{equation}
 and 
\begin{equation}
\left|N^{t-1}\left(A_{+}\right)\cup N^{t}\left(A_{-}\right)\right|\geq\left|N^{t}\left(A_{-}\right)\right|\geq\left|N^{t}\left(C_{-}\right)\right|.\label{eq:L5}
\end{equation}
Hence combining (\ref{eq:L2}), (\ref{eq:L3}), (\ref{eq:L4}) and
(\ref{eq:L5}) it follows that 
\begin{equation}
\left|N^{t}\left(A\right)\right|\ge\left|N^{t}\left(C\right)\right|.\label{eq:L6}
\end{equation}
However, since $\left|C\right|=\left|A\right|$ and $N^{t}\left(A\right)$
is minimal, the equality must hold in (\ref{eq:L6}), and in particular
it must hold in (\ref{eq:L4}) and (\ref{eq:L5}). Thus (\ref{eq:L4})
implies that we have $\left|N^{t}\left(A_{+}\right)\right|=\left|N^{t}\left(C_{+}\right)\right|$,
and since $C_{+}$ is the initial segment of the simplicial order,
it follows that $N^{t}\left(A_{+}\right)$ is minimal. 
Since the equality holds in (\ref{eq:L4}), we also obtain that $N^{t-1}\left(A_{-}\right)\cup N^{t}\left(A_{+}\right)=N^{t}\left(A_{+}\right)$.
Hence it follows that $N^{t-1}\left(A_{-}\right)\subseteq N^{t}\left(A_{+}\right)$.
Similarly the fact that the equality holds in (\ref{eq:L5}) implies that
$N^{t}\left(A_{-}\right)$ is minimal, and that $N^{t-1}\left(A_{+}\right)\subseteq N^{t}\left(A_{-}\right)$.
\end{proof}
\begin{prop}
\label{prop:Prop 11}Suppose that $A\subseteq Q_{n}$ is an extremal
set for which there exist $t\in Q_{n}$ and $r$ such that $B\left(t,r\right)\subseteq A\subseteq B\left(t,r+1\right)$
. Then $A$ is isomorphic to the initial segment of the simplicial
order. 
\end{prop}
\begin{proof}
The proof is by induction on $n$. When $n\leq2$ it is easy to verify
that the claim is true. Suppose that the claim holds for $n-1$, and
let $A\subseteq Q_{n}$ be an extremal set for which there exist $t\in Q_{n}$
and $r$ such that $B\left(t,r\right)\subseteq A\subseteq B\left(t,r+1\right)$.
By the symmetry of $Q_{n}$, we may assume that $t=\emptyset$. 

If $A$ satisfies $\left|A\right|=f_{r}$, then Proposition \ref{prop:Prop5}
implies that $A$ is an exact Hamming ball of radius $r$. In particular,
$A$ is isomorphic to the initial segment of the simplicial order,
and the claim follows in this case. Otherwise, by taking complements
if necessary, we may assume that $A$ satisfies $f_{r}<\left|A\right|\leq g_{r}$.
Indeed, this follows from the earlier observation that extremality
is preserved under taking complements. 

Since $f_{r}<\left|A\right|\leq g_{r}$, Lemma \ref{lem:Lemma6}  implies
that there exist $x,\,y$ and $z$ with $y\neq z$ satisfying $B\left(x,r\right)\subseteq A\subseteq B\left(y,r+1\right)\cap B\left(z,r+1\right)$,
$d\left(x,y\right)\leq1$ and $d\left(x,z\right)\leq1$ . We split
the rest of the proof into two cases based on whether $\left|A\right|=g_{r}$
or $f_{r}<\left|A\right|<g_{r}$.
\addtocounter{thm}{-11}
\begin{case}
The size of $A$ satisfies $\left|A\right|=g_{r}$.
\end{case}
\addtocounter{thm}{10}
Since $A\subseteq B\left(y,r+1\right)\cap B\left(z,r+1\right)$, it
follows that $B\left(y^{c},n-r-2\right)\cup B\left(z^{c},n-r-2\right)\subseteq A^{c}$.
Since $y^{c}\neq z^{c}$, Lemma \ref{lem:Lemma 7} implies that $\left|A^{c}\right|\geq g_{n-r-2}$.
But $\left|A^{c}\right|=2^{n}-g_{r}=g_{n-r-2}$, so the equality must
hold. In particular, this implies that we must have $A^{c}=B\left(y^{c},n-r-2\right)\cup B\left(z^{c},n-r-2\right)$,
and hence $A=B\left(y,r+1\right)\cap B\left(z,r+1\right)$. 

If $d\left(y,z\right)=1$, by applying a suitable isometry of $Q_{n}$
if necessary we may assume that $y=\emptyset$ and $z=\left\{ 1\right\} $.
Hence $A$ equals $B\left(\emptyset,r+1\right)\cap B\left(\left\{ 1\right\} ,r+1\right)=X^{\left(\leq r\right)}\cup\left(\left\{ 1\right\} +X_{1}^{\left(r\right)}\right)$,
which is the initial segment of the simplicial order of size $g_{r}$.
If $d\left(y,z\right)=2$, note that the only elements $s$ satisfying
the condition $A\subseteq B\left(s,r+1\right)$ are $s=y$ and $s=z$.
However, since $d\left(y,z\right)=2$ it follows that neither of $B\left(y,r\right)$
nor $B\left(z,r\right)$ is a subset of $A$. This contradicts the
existence of $t$ with $B\left(t,r\right)\subseteq A\subseteq B\left(t,r+1\right)$.
Hence $A$ has to be isomorphic to the initial segment of the simplicial
order. $\hfill\square$

\addtocounter{thm}{-10}
\begin{case}
The size of $A$ satisfies $f_{r}<\left|A\right|<g_{r}$.
\end{case}
\addtocounter{thm}{9}

Recall that $t=\emptyset$ is the element satisfying the condition
$B\left(t,r\right)\subseteq A\subseteq B\left(t,r+1\right)$, and
that $x,\,y$ and $z$ are elements with $y\neq z$ satisfying $B\left(x,r\right)\subseteq A\subseteq B\left(y,r+1\right)\cap B\left(z,r+1\right)$,
$d\left(x,y\right)\leq1$ and $d\left(x,z\right)\leq1$. If $x\neq\emptyset$,
then Lemma \ref{lem:Lemma 7} implies that $\left|A\right|\geq\left|B\left(\emptyset,r\right)\cup B\left(x,r\right)\right|\geq g_{r}$,
which contradicts the assumption $\left|A\right|<g_{r}$. Hence we
must have $x=\emptyset$. 

Since $y\neq z$, it follows that at least one of $y$ and $z$ does
not equal $\emptyset$, and hence we may assume that $y\ne\emptyset$.
Thus $d\left(y,x\right)\le1$ implies that $y=\left\{ i\right\} $
for some $i$, and by applying a suitable isometry if necessary we
may assume that $i=1$. Since $A\subseteq B\left(y,r+1\right)$, and
$B\left(\emptyset,r\right)\subseteq A\subseteq B\left(\emptyset,r+1\right)$,
it follows that $B\left(\emptyset,r\right)\subseteq A\subseteq B\left(\emptyset,r+1\right)\cap B\left(\left\{ 1\right\} ,r+1\right)$.
Thus there exists ${\cal B}\subseteq X_{1}^{\left(r\right)}$ for
which $A=X^{\left(\le r\right)}\cup\left(\left\{ 1\right\} +{\cal B}\right)$. 

Consider the $j$-sections in the direction $j=1$, and for simplicity
denote them by $A_{+}$ and $A_{-}$ for the rest of this proof. Since
$A=X^{\left(\le r\right)}\cup\left(\left\{ 1\right\} +{\cal B}\right)$,
it follows that $A_{+}=X_{1}^{\left(\leq r-1\right)}\cup{\cal B}$
and $A_{-}=X_{1}^{\left(\leq r\right)}$. Let $C_{+}$ and $C_{-}$
be the initial segments of the simplicial order of same sizes as $A_{+}$
and $A_{-}$ on ${\cal P}\left(X_{1}\right)$.

Since $f_{n-1,\,r-1}\leq\left|C_{+}\right|<\left|C_{-}\right|=f_{n-1,\,r}$,
it follows that $C_{-}\subseteq N\left(C_{+}\right)$. Since $A$
is an extremal set, it follows that $N^{t}\left(A\right)$ are minimal
for all $t>0$. Hence Lemma \ref{lem:Lemma 10} applied for all $t>0$
implies that $N^{t}\left(A_{+}\right)$ and $N^{t}\left(A_{-}\right)$
are minimal for all $t>0$. 

Let $D_{-}$ and $D_{+}$ be the initial segments of the simplicial
order of same sizes as $A_{-}^{c}$ and $A_{+}^{c}$ on ${\cal P}\left(X_{1}\right)$.
Since $f_{n-1,\,n-r-2}=\left|A_{-}^{c}\right|<\left|A_{+}^{c}\right|\leq f_{n-1,\,n-r-1}$,
it follows that $D_{+}^{c}\subseteq N\left(D_{-}^{c}\right)$. Hence
Lemma \ref{lem:Lemma 10} implies that $N^{t}\left(A_{+}^{c}\right)$
and $N^{t}\left(A_{-}^{c}\right)$ are minimal for all $t>0$, and
thus both $A_{+}$ and $A_{-}$ are extremal, although note that the
extremality of $A_{-}$ is also evident from the fact that $A_{-}=X_{1}^{\left(\leq r\right)}$. 

The extremality of $A_{+}$ implies that $X_{1}^{\left(\le r-1\right)}\cup{\cal B}$
is extremal as a subset of ${\cal P}\left(X_{1}\right)$. Since $B_{1}\left(\emptyset,r-1\right)\subseteq X_{1}^{\left(\le r-1\right)}\cup{\cal B}\subseteq B_{1}\left(\emptyset,r\right)$,
the inductive hypothesis implies that $X_{1}^{\left(\le r-1\right)}\cup{\cal B}$
is isomorphic to the initial segment of the simplicial order on ${\cal P}\left(X_{1}\right)$.
Recall that the initial segment of the simplicial order of size $\left|A_{+}\right|$
on ${\cal P}\left(X_{1}\right)$ is of the form $X_{1}^{\left(\leq r-1\right)}\cup{\cal C}$,
where ${\cal C}\subseteq X_{1}^{\left(r\right)}$ is the initial segment
of the lexicographic order on $X_{1}^{\left(r\right)}$ of size $\left|{\cal B}\right|$.
Thus Lemma \ref{lem:Lemma 9} implies that one can choose the isomorphism
to be of the form $\phi_{\sigma}$ for some $\sigma\in S_{n}$. Hence
we have $\phi_{\sigma}\left({\cal B}\right)={\cal C}$, which implies
that ${\cal B}$ is isomorphic to the initial segment of the lexicographic
order on $X_{1}^{\left(r\right)}$ as subsets of $X_{1}^{\left(r\right)}$.
Thus $\left\{ 1\right\} +{\cal B}$ is isomorphic to the initial segment
of the lexicographic order on $X^{\left(r+1\right)}$ as subsets of
$X^{\left(r+1\right)}$, and hence $A=X^{\left(\leq r\right)}\cup\left(\left\{ 1\right\} +{\cal B}\right)$
is isomorphic to the initial segment of the simplicial order, as required. 
\end{proof}
As usual, define the lower shadow of a set system ${\cal A}\subseteq X^{\left(r\right)}$
by $\partial^{-}{\cal A}=\left\{ b\setminus\left\{ i\right\} \,:\,b\in{\cal A},\,i\in b\right\} $,
and the iterated lower shadow by $\partial^{-t}{\cal A}=\partial^{-}\left(\partial^{-(t-1)}{\cal A}\right)$.
Similarly define the upper shadow of ${\cal A}\subseteq X^{\left(r\right)}$
by  $\partial^{+}{\cal A}=\left\{ b\cup\left\{ i\right\} \,:\,i\in X\setminus b,\,b\in{\cal A}\right\} $,
and the iterated upper shadow by $\partial^{+t}{\cal A}=\partial^{+}\left(\partial^{+(t-1)}{\cal A}\right)$.
Note that the upper shadow depends on the ground set, which will be
$X$ unless otherwise highlighted in the notation. For ${\cal A}\subseteq X^{\left(r\right)}$
define ${\cal \overline{A}}=\left\{ a^{c}\,:\,a\in{\cal A}\right\} $.
Note that $\left|{\cal \overline{A}}\right|=\left|{\cal A}\right|$
and ${\cal \overline{A}}\subseteq X^{\left(n-r\right)}$. It turns
out that the upper and lower shadows can be related to each others
via $\partial^{+}{\cal A}=\overline{\partial^{-}{\cal \overline{A}}}$. 

It is natural to ask that how should one choose a set system ${\cal A}\subseteq X^{\left(r\right)}$
of given size in order to minimise the size of the lower shadow of
${\cal A}$. Note that answering this question would also give an
answer for the same question concerning upper shadows by the earlier
observation. These questions are answered by the Kruskal-Katona theorem. 
\begin{thm}[Kruskal-Katona theorem \cite{key-15,key-17}]
\label{thm:Theorem 12 KK}
~\smallskip

1. Let ${\cal A}\subseteq X^{\left(r\right)}$, and let ${\cal B}\subseteq X^{\left(r\right)}$
be the initial segment of the colexicographic order with $\left|{\cal B}\right|=\left|{\cal A}\right|$.
Then $\left|\partial^{-}{\cal A}\right|\geq\left|\partial^{-}{\cal B}\right|$. 

2. Let ${\cal A}\subseteq X^{\left(r\right)}$, and let ${\cal C}\subseteq X^{\left(r\right)}$
be the initial segment of the lexicographic order with $\left|{\cal C}\right|=\left|{\cal A}\right|$.
Then $\left|\partial^{+}{\cal A}\right|\geq\left|\partial^{+}{\cal C}\right|$. $\hfill\square$
\end{thm}
Note that the lower shadow of an initial segment of the colexicographic
order on $X^{\left(r\right)}$ is an initial segment of the colexicographic
order on $X^{\left(r-1\right)}$, and similarly the upper shadow of
an initial segment of the lexicographic order on $X^{\left(r\right)}$
is an initial segment of the lexicographic order on $X^{\left(r+1\right)}$.
Hence we can strengthen the conclusion of Theorem \ref{thm:Theorem 12 KK}
by replacing $\left|\partial^{-}{\cal A}\right|\geq\left|\partial^{-}{\cal B}\right|$
with $\left|\partial^{-t}{\cal A}\right|\geq\left|\partial^{-t}{\cal B}\right|$
for all $t>0$, and by replacing $\left|\partial^{+}{\cal A}\right|\geq\left|\partial^{+}{\cal C}\right|$
with $\left|\partial^{+t}{\cal A}\right|\geq\left|\partial^{+t}{\cal C}\right|$
for all $t>0$. 

Let ${\cal A}\subseteq X^{\left(r\right)}$, let ${\cal B}$ be the
initial segment of the colexicographic order on $X^{\left(r\right)}$
of size $\left|{\cal A}\right|$ and let ${\cal C}$ be the initial
segment of the lexicographic order on $X^{\left(r\right)}$ of size
$\left|{\cal A}\right|$. We say that $\partial^{-t}{\cal A}$ are
\textit{minimal for all $t>0$} if $\left|\partial^{-t}{\cal A}\right|=\left|\partial^{-t}{\cal B}\right|$
for all $t>0$, and we say that $\partial^{+t}{\cal A}$ are \textit{minimal
for all $t>0$} if $\left|\partial^{+t}{\cal A}\right|=\left|\partial^{+t}{\cal C}\right|$
for all $t>0$. Since the Kruskal-Katona theorem implies that we always
have $\left|\partial^{-t}{\cal A}\right|\geq\left|\partial^{-t}{\cal B}\right|$
and $\left|\partial^{+t}{\cal A}\right|\geq\left|\partial^{+t}{\cal C}\right|$
for all $t>0$, in order to verify minimality it suffices to prove
that $\left|\partial^{-t}{\cal B}\right|\geq\left|\partial^{-t}{\cal A}\right|$
and $\left|\partial^{+t}{\cal C}\right|\geq\left|\partial^{+t}{\cal A}\right|$
for all $t>0$. 

Let ${\cal C}$ be an initial segment of the lexicographic order on
$X^{\left(r\right)}$, and let ${\cal B}=X^{\left(r\right)}\setminus{\cal C}$.
Then ${\cal B}$ is isomorphic to an initial segment of the colexicographic
order as subsets of $X^{\left(r\right)}$. Hence it follows that $\partial^{+t}{\cal C}$
are minimal for all $t>0$ and $\partial^{-t}{\cal B}$ are minimal
for all $t>0$. We now prove that if ${\cal A}\subseteq X^{\left(r\right)}$
for which $\partial^{+t}{\cal A}$ are minimal for all $t>0$ and
$\partial^{-t}\left(X^{\left(r\right)}\setminus{\cal A}\right)$ are
minimal for all $t>0$, then ${\cal A}$ and the initial segment of
the lexicographic order on $X^{\left(r\right)}$ of size $\left|{\cal A}\right|$
are isomorphic as subsets of $X^{\left(r\right)}$. 
\begin{cor}
\label{cor:Corollary 13}Let ${\cal A}\subseteq X^{(r)}$ and set
${\cal D}=X^{(r)}\setminus{\cal A}$. Suppose that $\partial^{+t}{\cal A}$
and $\partial^{-t}{\cal D}$ are minimal for all $t>0$. Then ${\cal A}$
and the initial segment of the lexicographic order of size $\left|{\cal A}\right|$
are isomorphic as subsets of $X^{\left(r\right)}$. 
\end{cor}
\begin{proof}
Let $A=X^{\left(\leq r-1\right)}\cup{\cal A}$, and note that $A^{c}=X^{\left(\geq r+1\right)}\cup{\cal D}$.
Let ${\cal C}\subseteq X^{\left(r\right)}$ be the initial segment
of the lexicographic order of size $\left|{\cal A}\right|$. Then
$C=X^{\left(\leq r-1\right)}\cup{\cal C}$ is the initial segment
of the simplicial order of size $\left|A\right|$. Define ${\cal B}=X^{\left(r\right)}\setminus{\cal C}$,
and note that $C^{c}=X^{\left(\geq r+1\right)}\cup{\cal B}$. 

Since $\partial^{+t}{\cal A}$ are minimal for all $t>0$, and ${\cal C}$
is the initial segment of the lexicographic order on $X^{\left(r\right)}$
with $\left|{\cal C}\right|=\left|{\cal A}\right|$, it follows that
$\left|\partial^{+t}{\cal A}\right|=\left|\partial^{+t}{\cal C}\right|$
for all $t>0$. Since $N^{t}\left(A\right)=X^{\left(\leq r+t-1\right)}\cup\partial^{+t}{\cal A}$
and $N^{t}\left(C\right)=X^{\left(\leq r+t-1\right)}\cup\partial^{+t}{\cal C}$,
it follows that 
\begin{equation}
\left|N^{t}\left(A\right)\right|=\left|X^{\left(\leq r+t-1\right)}\right|+\left|\partial^{+t}{\cal A}\right|=\left|X^{\left(\leq r+t-1\right)}\right|+\left|\partial^{+t}{\cal C}\right|=\left|N^{t}\left(C\right)\right|\label{eq:10'}
\end{equation}
 for all $t>0$. 

Since $\partial^{-t}{\cal D}$ are minimal for all $t>0$, and ${\cal B}$
is isomorphic to the initial segment of the colexicographic order
of size $\left|{\cal D}\right|$ as subsets of $X^{\left(r\right)}$,
it follows that $\left|\partial^{-t}{\cal D}\right|=\left|\partial^{-t}{\cal B}\right|$
for all $t>0$. Note that similarly we also have $N^{t}\left(A^{c}\right)=X^{\left(\geq r-t+1\right)}\cup\partial^{-t}{\cal D}$
and $N^{t}\left(C^{c}\right)=X^{\left(\geq r-t+1\right)}\cup\partial^{-t}{\cal B}$.
Hence it follows that 
\begin{equation}
\left|N^{t}\left(A^{c}\right)\right|=\left|X^{\left(\ge r-t+1\right)}\right|+\left|\partial^{-t}{\cal D}\right|=\left|X^{\left(\geq r-t+1\right)}\right|+\left|\partial^{-t}{\cal B}\right|=\left|N^{t}\left(C^{c}\right)\right|\label{eq:11'}
\end{equation}
for all $t>0$. Since (\ref{eq:10'}) and (\ref{eq:11'}) hold for
all $t>0$, it follows that $A$ is extremal. 

Since $B\left(\emptyset,r-1\right)\subseteq A\subseteq B\left(\emptyset,r\right)$,
Proposition \ref{prop:Prop 11} implies that $A$ is isomorphic to
the initial segment of the simplicial order, and hence $A$ is isomorphic
to $C$. Thus Lemma \ref{lem:Lemma 9} implies that there exists $\sigma\in S_{n}$
for which $\phi_{\sigma}\left(A\right)=C$. Since $\phi_{\sigma}$
maps the elements of $X^{\left(r\right)}$ to the elements of $X^{\left(r\right)}$,
it follows that $\phi_{\sigma}\left({\cal A}\right)={\cal C}$. Hence
${\cal A}$ and ${\cal C}$ are isomorphic as subsets of $X^{\left(r\right)}$,
as required. 
\end{proof}
For convenience, we recall the definition of the sets $A_{i}$ and
restate Theorem \ref{thm:Theorem 2}. Define the maps $\pi_{i}\,:\,X^{\left(r+1\right)}\rightarrow X_{i}^{\left(r+1\right)}\cup X_{i}^{\left(r\right)}$\textit{
}by $\pi_{i}\left(x\right)=x\setminus\left\{ i\right\} $ for all
$x\in X^{\left(r+1\right)}$, and for a set system ${\cal B}\subseteq X^{\left(r+1\right)}$
define $\pi_{i}\left({\cal B}\right)=\left\{ \pi_{i}\left(x\right)\,:\,x\in{\cal B}\right\} $.
Let $s$ be an integer of the form $s=f_{r}+k$ for some $0\leq k\leq{n-1 \choose r}$.
Let ${\cal A}$ be the initial segment of the lexicographic order
on $X^{\left(r+1\right)}$ of size $k$. Finally for each $i$ set
\[
A_{i}=X^{\left(\leq r\right)}\cup\left(\left\{ i\right\} +\pi_{i}\left({\cal A}\right)\right),
\]
and note that $\left|A_{i}\right|=s$ for all $i$. 
\addtocounter{thm}{-12}
\begin{thm}
Let $A\subseteq Q_{n}$ be a subset of size $s$,
where $s=f_{r}+k$ for some $r$ and $k\leq{n-1 \choose r}$. Let
$A_{1},\dots,A_{n}$ be the sets defined as above for these choices
of $r$ and $k$. Then $A$ is extremal if and only if $A$ is isomorphic
to some $A_{i}$. 
\end{thm}
\addtocounter{thm}{11}
\begin{proof}
Let $0\leq r\leq n$, $0\leq k\le{n-1 \choose r}$ and set $s=f_{r}+k$.
Let $A_{1},\dots,A_{n}$ be the subsets of $Q_{n}$ defined as above
for these choices of $n,\,r,\,k$ and $s$. We first prove that every
extremal set $A\subseteq Q_{n}$ of size $s$ is isomorphic to one
of $A_{i}$, and then that each $A_{i}$ is extremal. The second part
is much easier, but since it follows very quickly as a consequence
of the first part, we start with the harder direction. 

If $k=0$, then $\left|A\right|=f_{r}$. Hence Proposition \ref{prop:Prop5}
implies that $A$ and each $A_{i}$ is an exact Hamming ball of radius
$r$. Hence the result follows in this case. Thus we may assume that
$k>0$. 

Since $f_{r}<\left|A\right|\leq g_{r}$, Lemma \ref{lem:Lemma6}  implies
that there exist $x,\,y,\,z\in Q_{n}$ with $y\ne z$ satisfying $d\left(x,y\right)\leq1$,
$d\left(x,z\right)\leq1$ and $B\left(x,r\right)\subseteq A\subseteq B\left(y,r+1\right)\cap B\left(z,r+1\right)$.
Without loss of generality assume that $x=\emptyset$. If $y=\emptyset$
or $z=\emptyset$, then $B\left(\emptyset,r\right)\subseteq A\subseteq B\left(\emptyset,r+1\right)$
and hence Proposition \ref{prop:Prop 11} implies that $A$ is isomorphic
to the initial segment of the simplicial order. Hence $A$ is isomorphic
to $A_{1}$, as required. Thus we may assume that $y\neq\emptyset$
and $z\neq\emptyset$. Hence $y=\left\{ i\right\} $ and $z=\left\{ j\right\} $
for distinct elements $i,\,j\in X$. 

The condition $B\left(x,r\right)\subseteq A\subseteq B\left(y,r+1\right)\cap B\left(z,r+1\right)$
guaranteed by Lemma \ref{lem:Lemma6} can be rewritten as 
\[
X^{\left(\leq r\right)}\subseteq A\subseteq X^{\left(\leq r\right)}\cup\left(\left\{ i,j\right\} +\left(X_{i,j}^{\left(r-1\right)}\cup X_{i,j}^{\left(r\right)}\right)\right).
\]
Hence $A=X^{\left(\leq r\right)}\cup\left(\left\{ i,j\right\} +{\cal A}_{1}\right)\cup\left(\left\{ i,j\right\} +{\cal A}_{2}\right)$,
where ${\cal A}_{t}\subseteq X_{i,j}^{\left(r-2+t\right)}$ for both
$t\in\left\{ 1,2\right\} $. Define ${\cal F}\subseteq X^{\left(r+1\right)}$
by setting ${\cal F}=\left(\left\{ i,j\right\} +{\cal A}_{1}\right)\cup\left(\left\{ i\right\} +{\cal A}_{2}\right)$,
and $C\subseteq Q_{n}$ by setting $C=X^{\left(\leq r\right)}\cup{\cal F}$.
Also note that ${\cal F}$ can be written as ${\cal F}=\left\{ i\right\} +\left(\left(\left\{ j\right\} +{\cal A}_{1}\right)\cup{\cal A}_{2}\right)$. 

Suppose that $C$ is extremal. Then by Proposition \ref{prop:Prop 11}
$C$ is isomorphic to the initial segment of the simplicial order
of size $\left|A\right|$, and by Lemma \ref{lem:Lemma 9} this isomorphism
can be chosen to be of the form $\phi_{\sigma}$ for some $\sigma\in S_{n}$.
Hence $\sigma\left({\cal F}\right)$ is the initial segment of the
lexicographic order on $X^{\left(r+1\right)}$ of size $\left|{\cal F}\right|$. 

Suppose that $\sigma\left(i\right)\neq1$. Since ${\cal F}\subseteq\left\{ i\right\} +X_{i}^{\left(r\right)}$,
it follows that $\sigma\left({\cal F}\right)\subseteq\sigma\left(i\right)+X_{\sigma\left(i\right)}^{\left(r\right)}$,
and since $\sigma\left({\cal F}\right)$ is an initial segment on
$X^{\left(r+1\right)}$ with $\left|\sigma\left({\cal F}\right)\right|\leq{n-1 \choose r}$,
it also follows that $\sigma\left({\cal F}\right)\subseteq\left\{ 1\right\} +X_{1}^{\left(r\right)}$.
Hence every element of $\sigma\left({\cal F}\right)$ contains $\left\{ 1,\sigma\left(i\right)\right\} $
as a subset. Thus $\sigma'$ given by $\sigma'=\left(1\,\sigma\left(i\right)\right)\sigma$
also maps ${\cal F}$ to the initial segment of the lexicographic
order, and hence we may assume that $\sigma\left(i\right)=1$. For
$\sigma\in S_{n}$ satisfying $\sigma\left(i\right)=1$ and which
maps ${\cal F}$ to the initial segment of the lexicographic order,
it is easy to see that we have $\phi_{\sigma}\left(A\right)=A_{\sigma\left(j\right)}$.
Hence in order to prove that $A$ is isomorphic to one of $A_{1},\dots,A_{n}$,
it suffices to show that $C$ is extremal. 

Let $\partial_{i,j}^{+}$ be the upper shadow operator with respect
to the ground set $X_{i,j}$. Note that we have 
\[
\partial^{+t}{\cal F}=\left(\left\{ i,j\right\} +\left(\partial_{i,j}^{+t}{\cal A}_{1}\cup\partial_{i,j}^{+\left(t-1\right)}{\cal A}_{2}\right)\right)\cup\left(\left\{ i\right\} +\partial_{i,j}^{+t}{\cal A}_{2}\right)
\]
and hence it follows that

\begin{equation}
\left|N^{t}\left(C\right)\right|=\left|X^{\left(\leq r+t\right)}\right|+\left|\partial^{+t}{\cal F}\right|=\left|X^{\left(\leq r+t\right)}\right|+\left|\partial_{i,j}^{+t}{\cal A}_{1}\cup\partial_{i,j}^{+\left(t-1\right)}{\cal A}_{2}\right|+\left|\partial_{i,j}^{+t}{\cal A}_{2}\right|.\label{eq:14}
\end{equation}
On the other hand we have
\[
N^{t}\left(A\right)=X^{\left(\leq r+t\right)}\cup\left(\left\{ i,j\right\} +\left(\partial_{i,j}^{+t}{\cal A}_{1}\cup\partial_{i,j}^{+\left(t-1\right)}{\cal A}_{2}\right)\right)\cup\left(\left\{ i,j\right\} +\partial_{i,j}^{+t}{\cal A}_{2}\right)
\]
and hence it follows that
\begin{equation}
\left|N^{t}\left(A\right)\right|=\left|X^{\left(\leq r+t\right)}\right|+\left|\partial_{i,j}^{+t}{\cal A}_{1}\cup\partial_{i,j}^{+\left(t-1\right)}{\cal A}_{2}\right|+\left|\partial_{i,j}^{+t}{\cal A}_{2}\right|.\label{eq:15}
\end{equation}
Combining (\ref{eq:14}) with (\ref{eq:15}) we obtain that $\left|N^{t}\left(A\right)\right|=\left|N^{t}\left(C\right)\right|$
for all $t>0$. 

Let ${\cal B}_{1}=X_{i,j}^{\left(r-1\right)}\setminus{\cal A}_{1}$
and ${\cal B}_{2}=X_{i,j}^{\left(r\right)}\setminus{\cal A}_{2}$.
Let $S_{r}=\left\{ x\in X^{\left(r\right)}\,:\,\left\{ i,j\right\} \not\subseteq x\right\} $,
i.e.\ $S_{r}=X^{\left(r\right)}\setminus\left(\left\{ i,j\right\} +X_{i,j}^{\left(r-2\right)}\right)$.
It is easy to verify that 
\begin{equation}
A^{c}=X^{\left(\geq r+3\right)}\cup\left(\left\{ i,j\right\} +\left({\cal B}_{1}\cup{\cal B}_{2}\right)\right)\cup S_{r+1}\cup S_{r+2}.\label{eq:16}
\end{equation}

Note that $\partial^{-t}S_{m}=S_{m-t}$ for any $m$ and $t$. Also
note that 
\[
\left\{ i,j\right\} +\partial^{-t}{\cal B}_{m}\subseteq\partial^{-t}\left(\left\{ i,j\right\} +{\cal B}_{m}\right)\subseteq\left(\left\{ i,j\right\} +\partial^{-t}{\cal B}_{m}\right)\cup S_{r+m-t}
\]
holds for both $m\in\left\{ 1,2\right\} $. Hence (\ref{eq:16}) implies
that 
\begin{equation}
N^{t}\left(A^{c}\right)=X^{\left(\geq r-t+3\right)}\cup\left(\left\{ i,j\right\} +\left(\partial^{-t}{\cal B}_{2}\cup\partial^{-\left(t-1\right)}{\cal B}_{1}\right)\right)\cup\left(\left\{ i,j\right\} +\partial^{-t}{\cal B}_{1}\right)\cup S_{r-t+1}\cup S_{r-t+2}.\label{eq:17}
\end{equation}
Note that $X^{\left(\geq r-t+3\right)}$, $\left\{ i,j\right\} +\left(\partial^{-t}{\cal B}_{2}\cup\partial^{-\left(t-1\right)}{\cal B}_{1}\right)$,
$S_{r-t+2}$ , $\left\{ i,j\right\} +\partial^{-t}{\cal B}_{1}$ and
$S_{r-t+1}$ are pairwisely disjoint subsets of $Q_{n}$. Indeed,
the first one contains only elements of size $\geq r-t+3$, whereas
the next two contain only elements of size $r-t+2$ and the last two
contain only elements of size $r-t+1$. Also, $S_{r-t+2}$ and $S_{r-t+1}$
do not have any elements in $\left\{ i,j\right\} +{\cal P}\left(X_{i,j}\right)$,
whereas the second and fourth set system are contained in $\left\{ i,j\right\} +{\cal P}\left(X_{i,j}\right)$. 

Combining these observations together with (\ref{eq:17}), we obtain
that 
\begin{equation}
\left|N^{t}\left(A^{c}\right)\right|=\left|X^{\left(\geq r-t+3\right)}\right|+\left|\partial^{-t}{\cal B}_{2}\cup\partial^{-\left(t-1\right)}{\cal B}_{1}\right|+\left|\partial^{-t}{\cal B}_{1}\right|+\left|S_{r-t+1}\right|+\left|S_{r-t+2}\right|.\label{eq:18}
\end{equation}
Note that we have $\left|X_{i,j}^{\left(r-1\right)}\right|+\left|X_{i,j}^{\left(r\right)}\right|={n-2 \choose r-1}+{n-2 \choose r}={n-1 \choose r}$
for all $r$. Since $S_{r}=X^{\left(r\right)}\setminus\left(\left\{ i,j\right\} +X_{i,j}^{\left(r-2\right)}\right)$,
it follows that 
\begin{align}
 & \left|S_{r-t+1}\right|+\left|S_{r-t+2}\right|=\left|X^{\left(r-t+2\right)}\right|+\left|X^{\left(r-t+1\right)}\right|-\left|X_{i,j}^{\left(r-t-1\right)}\right|-\left|X_{i,j}^{\left(r-t\right)}\right|\label{eq:19}\\
 & =\left|X^{\left(r-t+2\right)}\right|+{n \choose r-t+1}-{n-1 \choose r-t}=\left|X^{\left(r-t+2\right)}\right|+{n-1 \choose r-t+1}.\nonumber 
\end{align}
Combining (\ref{eq:18}) and (\ref{eq:19}) together with the fact
that $\left|X_{i}^{\left(r-t+1\right)}\right|={n-1 \choose r-t+1}$,
it follows that 
\begin{equation}
\left|N^{t}\left(A^{c}\right)\right|=\left|X^{\left(\geq r-t+2\right)}\right|+\left|X_{i}^{\left(r-t+1\right)}\right|+\left|\partial^{-t}{\cal B}_{2}\cup\partial^{-\left(t-1\right)}{\cal B}_{1}\right|+\left|\partial^{-t}{\cal B}_{1}\right|.\label{eq:20}
\end{equation}

Let ${\cal G}=X^{\left(r+1\right)}\setminus{\cal F}$. Since ${\cal F}\subseteq X^{\left(r+1\right)}$
and $C=X^{\left(\leq r\right)}\cup{\cal F}$, it follows that $C^{c}=X^{\left(\geq r+2\right)}\cup{\cal G}$.
Hence we have 
\begin{equation}
\left|N^{t}\left(C^{c}\right)\right|=\left|X^{\left(\geq r-t+2\right)}\right|+\left|\partial^{-t}{\cal G}\right|\label{eq:21}
\end{equation}
for all $t>0$. We now find an expression for $\left|\partial^{-t}{\cal G}\right|$
in terms of ${\cal B}_{1}$ and ${\cal B}_{2}$. 

Recall that ${\cal F}$ can be written as ${\cal F}=\left\{ i\right\} +\left(\left(\left\{ j\right\} +{\cal A}_{1}\right)\cup{\cal A}_{2}\right)$.
Hence ${\cal G}=X_{i}^{\left(r+1\right)}\cup\left(\left\{ i,j\right\} +{\cal B}_{1}\right)\cup\left(\left\{ i\right\} +{\cal B}_{2}\right)$.
Hence it is easy to verify that 
\begin{equation}
\partial^{-t}{\cal G}=X_{i}^{\left(r-t+1\right)}\cup\left(\left\{ i,j\right\} +\partial^{-t}{\cal B}_{1}\right)\cup\left(\left\{ i\right\} +\left(\partial^{-\left(t-1\right)}{\cal B}_{1}\cup\partial^{-t}{\cal B}_{2}\right)\right).\label{eq:22}
\end{equation}
Note that the sets $X_{i}^{\left(r-t+1\right)}$, $\left\{ i,j\right\} +\partial^{-t}{\cal B}_{1}$
and $\left\{ i\right\} +\left(\partial^{-\left(t-1\right)}{\cal B}_{1}\cup\partial^{-t}{\cal B}_{2}\right)$
are pairwisely disjoint. Indeed, this follows by noting that the sets
in $X_{i}^{\left(r-t+1\right)}$ do not contain $i$, the sets in
$\left\{ i,j\right\} +\partial^{-t}{\cal B}_{1}$ contain both $i$
and $j$, and the sets in $\left\{ i\right\} +\left(\partial^{-\left(t-1\right)}{\cal B}_{1}\cup\partial^{-t}{\cal B}_{2}\right)$
contain $i$ but not $j$. Hence (\ref{eq:21}) and (\ref{eq:22})
imply that 
\begin{equation}
\left|N^{t}\left(C^{c}\right)\right|=\left|X^{\left(\geq r-t+2\right)}\right|+\left|X_{i}^{\left(r-t+1\right)}\right|+\left|\partial^{-t}{\cal B}_{1}\right|+\left|\partial^{-\left(t-1\right)}{\cal B}_{1}\cup\partial^{-t}{\cal B}_{2}\right|.\label{eq:23}
\end{equation}
In particular, (\ref{eq:20}) and (\ref{eq:23}) imply that $\left|N^{t}\left(C^{c}\right)\right|=\left|N^{t}\left(A^{c}\right)\right|$
for all $t>0$. Hence $C$ is extremal, and thus by earlier observations
it follows that $A$ is isomorphic to $A_{i}$ for some $i$. 

Conversely, suppose that $A=A_{i}$ for some $i$. Define ${\cal A}_{1}$,
${\cal A}_{2}$, ${\cal B}_{1}$, ${\cal B}_{2}$, ${\cal F}$ and
${\cal G}$ as before. Then by the construction of $A_{i}$'s, it
follows that ${\cal F}$ is isomorphic to the initial segment of the
lexicographic order as subsets of $X^{\left(r+1\right)}$. As before,
define $C=X^{\left(\leq r\right)}\cup{\cal F}$. Then $C$ is isomorphic
to the initial segment of the simplicial order, so $C$ is extremal.
Then (\ref{eq:14}) and (\ref{eq:15}) still give that $\left|N^{t}\left(A\right)\right|=\left|N^{t}\left(C\right)\right|$
for all $t>0$, and similarly (\ref{eq:20}) and (\ref{eq:23}) give
that $\left|N^{t}\left(A^{c}\right)\right|=\left|N^{t}\left(C^{c}\right)\right|$
for all $t>0$. Since $A$ and $C$ are subsets of the same size and
$C$ is extremal, it follows that $A$ is extremal as well, as required. 
\end{proof}
Note that if $s=f_{r}$ for some $r$, then Proposition \ref{prop:Prop5}
implies that if $A$ is an extremal set of size $s$, then $A$ is
isomorphic to the initial segment of the simplicial order. We now
show that for all other sizes $s$ there exists a non-trivial extremal
set of size $s$.
\begin{cor}
For all $n$ and for all $s\not\in\left\{ f_{0},\dots,f_{n}\right\} $
there exists an extremal set $A\subseteq Q_{n}$ of size $s$ which
is not isomorphic to the initial segment of the simplicial order. 
\end{cor}
\begin{proof}
If $s=g_{r}$ for some $r$, we may take $A=B_{r}$ where $B_{r}$
is the set defined in Section 2. Indeed, recall that in Section 2
we proved that $B_{r}$ is extremal and that $B_{r}$ is not isomorphic
to the initial segment of the simplicial order. Hence by taking complements
if necessary, we may assume that $f_{r}<s<g_{r}$ for some $r$. Let
$k=s-f_{r}$, and let ${\cal D}$ be the initial segment of the lexicographic
order of size $k$ on $X^{\left(r+1\right)}$. Let ${\cal D}_{i,\,+}$
and ${\cal D}_{i,\,-}$ be the $i$-sections of ${\cal D}$. Since
${\cal D}\neq\emptyset$, there exists $i$ for which ${\cal D}_{i,\,-}\neq\emptyset$. 

Consider the set $A_{i}$ for this particular choice of $i$. $A_{i}$
is defined by $A_{i}=X^{\left(\leq r\right)}\cup\left(\left\{ i\right\} +{\cal D}_{i,\,+}\right)\cup\left(\left\{ i\right\} +{\cal D}_{i,\,-}\right)$.
Since ${\cal D}_{i,\,-}\neq\emptyset$, it follows that $A_{i}\cap X^{\left(r+2\right)}\neq\emptyset$.
Since $s<g_{r}$, it follows that $B\left(\emptyset,r\right)$ is
the unique exact Hamming ball of radius $r$ contained in $A_{i}$.
Indeed, if $B\left(x,r\right)\subseteq A$ for some $x\neq\emptyset$,
then Lemma \ref{lem:Lemma 7} implies that $\left|A\right|\geq\left|B\left(\emptyset,r\right)\cup B\left(x,r\right)\right|\geq g_{r}$,
which contradicts $\left|A\right|<g_{r}$. Thus at least one of the
conditions $B\left(x,r\right)\subseteq A_{i}$ and $A_{i}\subseteq B\left(x,r+1\right)$
is violated for any $x\in Q_{n}$, since the second one is violated
for $x=\emptyset$ as $A_{i}\cap X^{\left(r+2\right)}\neq\emptyset$.
In particular, it follows that $A_{i}$ is not a Hamming ball. Hence
$A_{i}$ is an extremal set which is not isomorphic to the initial
segment of the simplicial order, as initial segments of the simplicial
order are always Hamming balls. 
\end{proof}
It is natural to ask that when are $A_{i}$ and $A_{j}$ isomorphic
as subsets of $Q_{n}$. Let ${\cal D}$ be the initial segment from
which $A_{i}$'s are obtained. If $\sigma=\left(ij\right)\in S_{n}$
satisfies $\phi_{\sigma}\left({\cal D}\right)={\cal D}$, then $A_{i}$
and $A_{j}$ are certainly isomorphic, where $\phi_{\sigma}$ is given
as before. The aim of the following lemma is to prove that this is
the only way the isomorphism can occur.
\begin{lem}
\label{lem:Lemma15}Let $0\leq r\leq n$, $0\leq k\leq{n-1 \choose r}$,
$s=f_{r}+k$, and let $A_{i}$ be the sets defined before Theorem
\ref{thm:Theorem 2} for these choices of $n$ and $s$. Let ${\cal D}$
be the initial segment of the lexicographic order of size $k$ on
$X^{\left(r+1\right)}$. Then $A_{i}$ and $A_{j}$ are isomorphic
if and only if $\phi_{\sigma}\left({\cal D}\right)={\cal D}$ for
$\sigma=\left(ij\right)$.
\end{lem}
We say that ${\cal D}\subseteq X^{\left(r\right)}$ is \textit{left
compressed} if for all $i<j$ and $a\in{\cal D}$ we have 
\[
\left(j\in a\text{ and }i\not\in a\right)\Rightarrow\left(a\setminus\left\{ j\right\} \right)\cup\left\{ i\right\} \in{\cal D}.
\]
Note that if ${\cal D}$ is an initial segment of the lexicographic
or colexicographic order, then ${\cal D}$ is left compressed. 
\begin{proof}
If $s=f_{r}$, then ${\cal D}=\emptyset$ and $A_{i}=B\left(\emptyset,r\right)$
for all $i$. Hence each $A_{i}$ and $A_{j}$ are isomorphic, and
$\phi_{\sigma}\left({\cal D}\right)={\cal D}$ for all $\sigma=\left(ij\right)$
and therefore the claim follows. If $s=g_{r}$, then ${\cal D}=\left\{ 1\right\} +X_{1}^{\left(r\right)}$.
Hence $\phi_{\sigma}\left({\cal D}\right)={\cal D}$ for $\sigma=\left(ij\right)$
if and only if $1\not\in\left\{ i,j\right\} $. On the other hand,
note that $A_{1}$ is isomorphic to the initial segment of the simplicial
order, and for $j\geq2$ each $A_{j}$ is isomorphic to the set $B_{r}$
defined in Section 2. Hence the claim follows if $s=f_{r}$ or $s=g_{r}$.
Thus we may assume that $f_{r}<s<g_{r}$. 

If $\phi_{\sigma}\left({\cal D}\right)={\cal D}$ for $\sigma=\left(ij\right)$,
it follows that $\phi_{\sigma}\left({\cal D}_{i,\,-}\right)={\cal D}_{j,\,-}$
and $\phi_{\sigma}\left({\cal D}_{i,\,+}\right)={\cal D}_{j,\,+}$.
Hence $\phi_{\sigma}\left(A_{i}\right)=A_{j}$, so $A_{i}$ and $A_{j}$
are isomorphic. 

Suppose that $A_{i}$ and $A_{j}$ are isomorphic and $i<j$. Since
$s<g_{r}$, Lemma \ref{lem:Lemma 7} implies that $B\left(\emptyset,r\right)$
is the unique exact Hamming ball of radius $r$ contained in each
$A_{t}$. Hence if $g$ is an isometry mapping $A_{i}$ to $A_{j}$,
then $g\left(\emptyset\right)=\emptyset$. Thus $g=\phi_{\sigma}$
for some $\sigma\in S_{n}$ by Lemma \ref{lem:Lemma 8}. Hence it
follows that $A_{i}\cap X^{\left(m\right)}$ is mapped to $A_{j}\cap X^{\left(m\right)}$
for all $m$, and in particular we have $\left|{\cal D}_{i,\,-}\right|=\left|{\cal D}_{j,\,-}\right|$
and $\left|{\cal D}_{i,\,+}\right|=\left|{\cal D}_{j,\,+}\right|$.
In order to prove that $\phi_{\left(ij\right)}\left({\cal D}\right)={\cal D}$,
we need to prove the following two claims

\begin{enumerate}
\item For all $a\in{\cal D}$ with $i\not\in a$ and $j\in a$ we have
$\left(a\setminus\left\{ j\right\} \right)\cup\left\{ i\right\} \in{\cal D}$.
\item For all $a\in{\cal D}$ with $i\in a$ and  $j\not\in a$ we have
$\left(a\setminus\left\{ i\right\} \right)\cup\left\{ j\right\} \in{\cal D}$.
\end{enumerate}

Since $i<j$ and ${\cal D}$ is an initial segment of the lexicographic
order, the first claim follows immediately from the fact that ${\cal D}$
is left-compressed. 

Define ${\cal T}_{i}=\left\{ a\in{\cal D}\,:\,i\in a,\,j\not\in a\right\} $,
${\cal T}_{j}=\left\{ a\in{\cal D}\,:\,j\in a,\,i\not\in a\right\} $
and ${\cal T}_{i,j}=\left\{ a\in{\cal D}\,:\,i\in a,\,j\in a\right\} $.
Note that ${\cal D}_{i,\,+}={\cal {\cal T}}_{i}\cup{\cal T}_{i,\,j}$
and ${\cal D}_{j,\,+}={\cal {\cal T}}_{j}\cup{\cal T}_{i,\,j}$. Since
${\cal T}_{i}$, ${\cal T}_{j}$ and ${\cal T}_{i,j}$ are disjoint
sets, the condition $\left|{\cal D}_{i,\,+}\right|=\left|{\cal D}_{j,\,+}\right|$
is equivalent to $\left|{\cal T}_{i}\right|=\left|{\cal T}_{j}\right|$. 

Since the first claim is true, it follows that the map $a\rightarrow\left(a\setminus\left\{ j\right\} \right)\cup\left\{ i\right\} $
is an injection from ${\cal T}_{j}$ to ${\cal T}_{i}$. Since $\left|{\cal T}_{i}\right|=\left|{\cal T}_{j}\right|$,
it follows that this map is bijection, and the inverse is given by
$b\rightarrow\left(b\setminus\left\{ i\right\} \right)\cup\left\{ j\right\} $.
Hence the second claim must be true as well, and thus we have $\phi_{\left(ij\right)}\left({\cal D}\right)={\cal D}$. 
\end{proof}
It is natural to ask whether for all $t$ there exist $n$ and $s$
for which there are $t$ pairwisely non-isomorphic extremal sets $B_{1},\dots,B_{t}$
of size $s$ on $Q_{n}$. One can use Lemma \ref{lem:Lemma15} to
conclude that this turns out to be true. Indeed, for given $t$ it
suffices to find $n$, $r$ and an initial segment ${\cal A}$ of
the lexicographic order on $X^{\left(r\right)}$ with $\left|{\cal A}\right|\leq{n-1 \choose r-1}$
for which the sizes of ${\cal A}_{j_{i},\,+}$ are distinct for some
$j_{1},\dots,j_{t}\in X$. 

Recall that if ${\cal A}$ is isomorphic to an initial segment of
the colexicographic order, then $X^{\left(r\right)}\setminus{\cal A}$
is isomorphic to an initial segment of the lexicographic order. Also
note that $\left|{\cal A}_{j,\,+}\right|+\left|\left(X^{\left(r\right)}\setminus{\cal A}\right)_{j,\,+}\right|={n-1 \choose r-1}$
for all $j\in X$. Hence it suffices to find an initial segment ${\cal A}$
of the colexicographic order on $X^{\left(r\right)}$ with $\left|{\cal A}\right|\geq{n \choose r}-{n-1 \choose r-1}={n-1 \choose r}$
for which the sizes of ${\cal A}_{j_{i},\,+}$ are distinct for some
$j_{1},\dots,j_{t}\in X$. 

For all $t$, we construct inductively ${\cal A}_{t}$ which is a
proper subset of $\left\{ 1,\dots,2t\right\} ^{\left(t\right)}$ satisfying
$\left|\left({\cal A}_{t}\right)_{2i,\,+}\right|<\left|\left({\cal A}_{t}\right)_{2j,\,+}\right|$
for all $1\leq j<i\leq t$ and $\left|{\cal A}_{t}\right|\geq{2t-1 \choose t}$.
As the base case $t=1$ we can certainly take ${\cal A}_{1}=\left\{ \left\{ 1\right\} \right\} $.

Now suppose that $t\geq2$ and that the claim holds for $t-1$, and
consider ${\cal A}_{t+1}=\left\{ 1,\dots,2t+1\right\} ^{\left(t+1\right)}\cup\left({\cal A}_{t}+\left\{ 2t+2\right\} \right)$.
We claim that this set system has the required properties. First of
all, ${\cal A}_{t+1}$ is certainly an initial segment of the colexicographic
order on $\left\{ 1,\dots,2t+2\right\} ^{\left(t+1\right)}$ as ${\cal A}_{t}$
is an initial segment of the colexicographic order on $\left\{ 1,\dots,2t\right\} ^{\left(t\right)}$,
and $\left|{\cal A}_{t+1}\right|\geq{2t+1 \choose t+1}$ as required. 

Note that for all $1\leq i\leq t$ we have $\left({\cal A}_{t+1}\right)_{2i,\,+}={2t \choose t}+\left|\left({\cal A}_{t}\right)_{2i,\,+}\right|\geq{2t \choose t}$.
Hence the inductive hypothesis implies that $\left|\left({\cal A}_{t+1}\right)_{2i,\,+}\right|<\left|\left({\cal A}_{t+1}\right)_{2j,\,+}\right|$
for all $1\leq j<i\leq t$. On the other hand, note that $\left|\left({\cal A}_{t+1}\right)_{2t+2,\,+}\right|=\left|{\cal A}_{t}\right|<{2t \choose t}$,
as ${\cal A}_{t}$ is a proper subset of $\left\{ 1,\dots,2t\right\} ^{\left(t\right)}$
by induction. Hence $\left|\left({\cal A}_{t+1}\right)_{2t+2,\,+}\right|<\left|\left({\cal A}_{t+1}\right)_{2i,\,+}\right|$
holds for all $1\leq i\leq t+1$, which completes the proof of the
inductive step. Hence for $n=2t$ there exists a suitable $s$ with
$f_{t-1}\leq s\leq f_{t-1}+{2t-1 \choose t-1}$ for which the sets
$A_{2},A_{4},\dots,A_{2t}$ are pairwisely non-isomorphic. 

\section{The weak version}

Suppose we weaken the notion of extremality so that we only require
$N\left(A\right)$ and $N\left(A^{c}\right)$ to have minimal size
among the sets $A\subseteq Q_{n}$ of given size. The aim of this
section is to prove that this weaker condition actually implies that
$\left|N^{t}\left(B\right)\right|\geq\left|N^{t}\left(A\right)\right|$
and $\left|N^{t}\left(B^{c}\right)\right|\geq\left|N^{t}\left(A^{c}\right)\right|$
for all $B\subseteq Q_{n}$ with $\left|B\right|=\left|A\right|$
and for all $t>0$, i.e.\ that $A$ is extremal. From now on, we
say that $A$ is \textit{weakly extremal }if for all $B\subseteq Q_{n}$
with $\left|B\right|=\left|A\right|$ we have $\left|N\left(B\right)\right|\geq\left|N\left(A\right)\right|$
and $\left|N\left(B^{c}\right)\right|\geq\left|N\left(A^{c}\right)\right|$. 

F\"{u}redi and Griggs \cite{key-12} proved that if ${\cal D}\subseteq X^{(r)}$
is a set system for which the size of $\partial^{-}{\cal D}$ is minimal
among the subsets of $X^{\left(r\right)}$ of given size, then $\partial^{-t}{\cal D}$
are minimal for all $t>0$. Our aim is to prove that a similar conclusion
holds for subsets of $Q_{n}$ as well, that is if $A\subseteq Q_{n}$
is a subset for which the size of $N\left(A\right)$ is minimal among
the subsets of $Q_{n}$ of a given size, then the size of $N^{t}\left(A\right)$
is also minimal for all $t>0$ among the subsets of $Q_{n}$ of that
given size. This result immediately implies that the notions of weak
extremality and extremality coincide. Hence the main task is to prove
the following theorem. 

\addtocounter{thm}{-12}
\begin{thm}
Let $A\subseteq Q_{n}$ be a subset for which
every $B\subseteq Q_{n}$ with $\left|B\right|=\left|A\right|$ satisfies
$\left|N\left(B\right)\right|\geq\left|N\left(A\right)\right|$. Then
for all $t>0$ and $B\subseteq Q_{n}$ with $\left|B\right|=\left|A\right|$
we also have $\left|N^{t}\left(B\right)\right|\geq\left|N^{t}\left(A\right)\right|$.
\end{thm}
\addtocounter{thm}{11}

Let $A\subseteq Q_{n}$ be a set satisfying the conditions of Theorem
\ref{thm:Theorem 4}, and let $A_{i,\,+}$ and $A_{i,\,-}$ be its
$i$-sections in some direction $i$. A natural way to prove Theorem
\ref{thm:Theorem 4} is to apply an inductive argument on $A_{i,\,+}$
and $A_{i,\,-}$. Let $C_{i,\,+}$ and $C_{i,\,-}$ be the initial
segments of the simplicial order on ${\cal P}\left(X_{i}\right)$
of sizes $\left|A_{i,\,+}\right|$ and $\left|A_{i,\,-}\right|$ respectively.
Without loss of generality we may assume that $\left|A_{i,\,+}\right|\leq\left|A_{i,\,-}\right|$.
If the sets $C_{i,\,+}$ and $C_{i,\,-}$ satisfy the conditions of
Lemma \ref{lem:Lemma 10} when $t=1$, i.e.\ $C_{i,\,-}\subseteq N\left(C_{i,\,+}\right)$
and $N\left(A\right)$ is minimal, then the Lemma 
implies that both $N\left(A_{i,\,+}\right)$ and $N\left(A_{i,\,-}\right)$
are minimal, and hence one can use induction on $n$ to prove Theorem
\ref{thm:Theorem 4}. However, we cannot apply this argument if the
conditions of Lemma \ref{lem:Lemma 10} are not satisfied when $t=1$, i.e.\ if
$N\left(C_{i,\,+}\right)$ is a proper subset of $C_{i,\,-}$ . It
is easy to verify that this may happen even when $A$ is an initial
segment of the simplicial order. For example, if $A=B\left(\emptyset,r\right)\cup\left\{ \left\{ 1,\dots,r+1\right\} \right\} $,
then $N\left(A_{i,\,+}\right)$ is a proper subset of $A_{i,\,-}$
for every $i>r+1$. 

Note that in the previous example there were some directions $i$
for which the sets $C_{i,\,-}$ and $C_{i,\,+}$ satisfy the conditions
of Lemma \ref{lem:Lemma 10}, for example any $i\leq r+1$ would work.
Conveniently, it turns out that the same holds in general: for any
set $A\subseteq Q_{n}$ satisfying the conditions of Theorem \ref{thm:Theorem 4}
there is always a direction $i$ for which both $C_{i,\,+}$ and $C_{i,\,-}$
satisfy the conditions of Lemma \ref{lem:Lemma 10}\@.

We start with some preliminary results, which are Lemmas \ref{lem:Lemma 16},
\ref{lem:Lemma 17} and \ref{lem:Lemma 18}. In a sense, Lemmas \ref{lem:Lemma 17}
and \ref{lem:Lemma 18} are just necessary tools for subsequent results,
and the proofs are mostly calculational. These statements are in flavor
similar to Harper's theorem and the Kruskal-Katona theorem, but it
seems that there is no straightforward way of deducing them directly
from Harper's theorem or the Kruskal-Katona theorem. 

Armed with those, we move on to Lemma \ref{lem:Lemma 19} in which
we prove that there exists a direction $i$ for which $C_{i,\,+}$
and $C_{i,\,-}$ satisfy the conditions of Lemma \ref{lem:Lemma 10}.
Given Lemma \ref{lem:Lemma 19}, the proof of Theorem \ref{thm:Theorem 4}
follows almost immediately, and finally we can deduce that weak extremality
implies extremality in Theorem \ref{thm:Theorem 3}.
\begin{lem}
\label{lem:Lemma 16}Let ${\cal A},\,{\cal B}\subseteq X^{\left(r\right)}$
be initial segments of the colexicographic order with $\left|{\cal A}\right|+\left|{\cal B}\right|\le{n \choose r}$,
and recall that ${\cal A}_{+}=\left\{ a\setminus\left\{ 1\right\} \,:\,1\in a,\,a\in{\cal A}\right\} $.
Let ${\cal C}\subseteq X^{\left(r\right)}$ be the initial segment
of the colexicographic order of size $\left|{\cal A}\right|+\left|{\cal B}\right|$
on $X^{\left(r\right)}$. Then $\left|{\cal A}_{+}\right|+\left|{\cal B}_{+}\right|\geq\left|{\cal C}_{+}\right|$. 
\end{lem}
\begin{proof}
Let ${\cal I}$ and ${\cal J}$ be initial segments of the colexicographic
order on $X^{\left(r\right)}$ with $\left|{\cal I}\right|\geq\left|{\cal J}\right|$,
and let ${\cal K}$ be the initial segment of the colexicographic
order on $X^{\left(r\right)}$ of size $\left|{\cal I}\setminus{\cal J}\right|$.
Note that since initial segments are nested, it follows that $\left|{\cal K}\right|+\left|{\cal J}\right|=\left|{\cal I}\right|$.
Lemma 9 in \cite{key-19} states that we have $\left|\left({\cal I\setminus{\cal J}}\right)_{+}\right|\leq\left|{\cal K}_{+}\right|$. 

We apply this Lemma for ${\cal I}={\cal C}$ and ${\cal J}={\cal A}$.
Hence we must have $\left|{\cal K}\right|=\left|{\cal B}\right|$,
and since the initial segments are unique it follows that ${\cal K}={\cal B}$.
Since ${\cal A}\subseteq{\cal C}$, it follows directly that $\left|{\cal C}_{+}\right|=\left|{\cal \left(C\setminus{\cal A}\right)}_{+}\right|+\left|{\cal A}_{+}\right|$.
Since the Lemma implies that $\left|\left({\cal C}\setminus{\cal A}\right)_{+}\right|\leq\left|{\cal B}_{+}\right|$,
it follows that $\left|{\cal A}_{+}\right|+\left|{\cal B}_{+}\right|\geq\left|{\cal C}_{+}\right|$.
\end{proof}
\begin{lem}
\label{lem:Lemma 17}1. Let ${\cal A},\,{\cal B}\subseteq X^{\left(r\right)}$
be non-empty initial segments of the lexicographic order with $\left|{\cal A}\right|+\left|{\cal B}\right|\leq{n \choose r}$,
and let ${\cal C}$ be the initial segment of the lexicographic order
of size $\left|{\cal A}\right|+\left|{\cal B}\right|$ on $X^{\left(r\right)}$.
Then $\left|\partial^{+}{\cal A}\right|+\left|\partial^{+}{\cal B}\right|>\left|\partial^{+}{\cal C}\right|$. 

2. Let ${\cal A},\,{\cal B}\subseteq X^{\left(r\right)}$ be non-empty
initial segments of the colexicographic order with $\left|{\cal A}\right|+\left|{\cal B}\right|\leq{n \choose r}$,
and let ${\cal C}$ be the initial segment of the colexicographic
order of size $\left|{\cal A}\right|+\left|{\cal B}\right|$ on $X^{\left(r\right)}$.
Then $\left|\partial^{-}{\cal A}\right|+\left|\partial^{-}{\cal B}\right|>\left|\partial^{-}{\cal C}\right|$. 
\end{lem}
Note that the Kruskal-Katona theorem certainly implies that $\left|\partial^{+}{\cal A}\right|+\left|\partial^{+}{\cal B}\right|\ge\left|\partial^{+}{\cal C}\right|$
and $\left|\partial^{-}{\cal A}\right|+\left|\partial^{-}{\cal B}\right|\ge\left|\partial^{-}{\cal C}\right|$.
However, there does not seem to be a straightforward way to directly
conclude the strict inequalities from the Kruskal-Katona theorem.
\begin{proof}
Let ${\cal A}$ be an initial segment of the lexicographic order and
let $\overline{{\cal A}}=\left\{ a^{c}\,:\,a\in{\cal A}\right\} $.
Recall that $\left|{\cal A}\right|=\left|{\cal \overline{A}}\right|$,
$\partial^{-}\overline{{\cal A}}=\overline{\partial^{+}{\cal A}}$
and that $\overline{{\cal A}}$ is isomorphic to an initial segment
of the colexicographic order. Hence it follows that these two claims
are equivalent, so it suffices to prove the second one. 

The proof is by induction on $r$. When $r=1$, $\left|\partial^{-}{\cal A}\right|=\left|\partial^{-}{\cal B}\right|=\left|\partial^{-}{\cal C}\right|=1$
as each of ${\cal A}$, ${\cal B}$ and ${\cal C}$ is non-empty,
and hence the claim follows trivially. Hence we may assume that $r>1$,
and that the claim holds for $r-1$. 

Recall that ${\cal A}_{+}=\left\{ a\setminus\left\{ 1\right\} \,:\,1\in a,\,a\in{\cal A}\right\} $
and ${\cal A}_{-}=\left\{ a\,:\,1\not\in a,\,a\in{\cal A}\right\} $.
Since $1$ is the smallest element in $X$ and ${\cal A}$ is an initial
segment of the colexicographic order and hence left compressed, it
follows that $\partial^{-}{\cal A}_{-}^{.}\subseteq{\cal A}_{+}$.
Since ${\cal A}=\left(\left\{ 1\right\} +{\cal A}_{+}\right)\cup{\cal A}_{-}$,
it follows that $\partial^{-}{\cal A}=\left(\left\{ 1\right\} +\partial^{-}{\cal A}_{+}\right)\cup\left({\cal A}_{+}\cup\partial^{-}{\cal A}_{-}\right)=\left(\left\{ 1\right\} +\partial^{-}{\cal A}_{+}\right)\cup{\cal A}_{+}$,
and hence 
\begin{equation}
\left|{\cal A}\right|=\left|{\cal A}_{+}\right|+\left|\partial^{-}{\cal A}_{+}\right|.\label{eq:24}
\end{equation}

Define ${\cal B}_{+}$ and ${\cal C}_{+}$ similarly. Clearly (\ref{eq:24})
holds for them as well. Since ${\cal A}$, ${\cal B}$ and ${\cal C}$
are non-empty initial segments of the colexicographic order, it follows
that ${\cal A}_{+}$, ${\cal B}_{+}$ and ${\cal C}_{+}$ are non-empty
initial segments of the colexicographic order on  $X_{1}^{\left(r-1\right)}$.
Lemma \ref{lem:Lemma 16} also implies that we have 
\begin{equation}
\left|{\cal A}_{+}\right|+\left|{\cal B}_{+}\right|\ge\left|{\cal C}_{+}\right|.\label{eq:25}
\end{equation}
Hence the inductive hypothesis implies that

\begin{equation}
\left|\partial^{-}{\cal A}_{+}\right|+\left|\partial^{-}{\cal B}_{+}\right|>\left|\partial^{-}{\cal C}_{+}\right|.\label{eq:26}
\end{equation}
Combining (\ref{eq:24}), (\ref{eq:25}) and (\ref{eq:26}), we obtain
that 
\begin{align*}
 & \left|\partial^{-}{\cal A}\right|+\left|\partial^{-}{\cal B}\right|=\left|\partial^{-}{\cal A}_{+}\right|+\left|\partial^{-}{\cal B}_{+}\right|+\left|{\cal A}_{+}\right|+\left|{\cal B}_{+}\right|\\
 & \geq\left|\partial^{-}{\cal A}_{+}\right|+\left|\partial^{-}{\cal B}_{+}\right|+\left|{\cal C}_{+}\right|>\left|\partial^{-}{\cal C}_{+}\right|+\left|{\cal C}_{+}\right|=\left|\partial^{-}{\cal C}\right|,
\end{align*}
which completes the proof of Lemma \ref{lem:Lemma 17}. 
\end{proof}
\begin{lem}
\label{lem:Lemma 18}Let $A\subseteq Q_{n}$ be a set such that for
all $B\subseteq Q_{n}$ with $\left|B\right|=\left|A\right|$ we have
$\left|N\left(B\right)\right|\geq\left|N\left(A\right)\right|$. Let
$D_{+}$ be the initial segment of the simplicial order on $Q_{n-1}$
of the largest size for which the initial segment $D_{-}$ of the
simplicial order on $Q_{n-1}$ of size $\left|A\right|-\left|D_{+}\right|$
satisfies $N\left(D_{+}\right)\subseteq D_{-}$. Then for any direction
$i$ we have $\left|A_{i,\,+}\right|\geq\left|D_{+}\right|$. Furthermore,
if $\left|N\left(A_{i,\,+}\right)\right|\leq\left|A_{i,\,-}\right|$
then we must have $\left|A_{i,\,+}\right|=\left|D_{+}\right|$. 
\end{lem}
\begin{proof}
Denote the $i$-sections by $A_{+}$ and $A_{-}$. Let $D_{+}$ and
$D_{-}$ be defined as in the statement, and let $D=D_{-}\cup\left(\left\{ n\right\} +D_{+}\right)$. 

Suppose that $\left|A_{+}\right|<\left|D_{+}\right|$, and hence we
also have $\left|A_{-}\right|>\left|D_{-}\right|$. Since $D_{-}$
is an initial segment of the simplicial order, Harper's theorem implies
that $\left|N\left(A_{-}\right)\right|\geq\left|N\left(D_{-}\right)\right|$.
Since $N\left(D_{+}\right)\subseteq D_{-}$, and hence $D_{+}\subseteq N^{2}\left(D_{+}\right)\subseteq N\left(D_{-}\right)$,
(\ref{eq:L2}) implies that 
\[
\left|N\left(D\right)\right|=\left|N\left(D_{-}\right)\right|+\left|D_{-}\right|.
\]
On the other hand, (\ref{eq:L2}) implies that 
\[
\left|N\left(A\right)\right|\geq\left|N\left(A_{-}\right)\right|+\left|A_{-}\right|>\left|N\left(D_{-}\right)\right|+\left|D_{-}\right|=\left|N\left(D\right)\right|,
\]
as $\left|A_{-}\right|>\left|D_{-}\right|$ and $\left|N\left(A_{-}\right)\right|\geq\left|N\left(D_{-}\right)\right|$.
This contradicts the fact that $\left|N\left(D\right)\right|\geq\left|N\left(A\right)\right|$.
Hence we must have $\left|A_{+}\right|\geq\left|D_{+}\right|$. 

In order to prove the second part, suppose that $\left|N\left(A_{+}\right)\right|\leq\left|A_{-}\right|$.
Let $C_{+}$ and $C_{-}$ be the initial segments of the simplicial
order on ${\cal P}\left(X_{i}\right)$ of sizes $\left|A_{+}\right|$
and $\left|A_{-}\right|$. Harper's theorem implies that $\left|N\left(A_{+}\right)\right|\geq\left|N\left(C_{+}\right)\right|$,
and since $\left|A_{-}\right|=\left|C_{-}\right|$ it follows that
$\left|N\left(C_{+}\right)\right|\leq\left|C_{-}\right|$. Since $C_{+}$
and $C_{-}$ are initial segments of the simplicial order on ${\cal P}\left(X_{i}\right)$,
it follows that $N\left(C_{+}\right)\subseteq C_{-}$, so $C_{+}$
and $C_{-}$ satisfies the conditions required from $D_{+}$ and $D_{-}$.
Thus by the maximality condition it follows that $\left|C_{+}\right|\leq\left|D_{+}\right|$.
By the first part, we know that $\left|A_{+}\right|\geq\left|D_{+}\right|$.
Since $\left|A_{+}\right|=\left|C_{+}\right|$, it follows that these
two inequalities imply that $\left|A_{+}\right|=\left|D_{+}\right|$. 
\end{proof}
\begin{lem}
\label{lem:Lemma 19}Let $A\subseteq Q_{n}$ be a set for which the
condition $\left|A\right|\neq f_{r}$ holds for all $r$, and such
that for all $B\subseteq Q_{n}$ with $\left|B\right|=\left|A\right|$
we have $\left|N\left(B\right)\right|\geq\left|N\left(A\right)\right|$.
Let $D_{+}$ and $D_{-}$ be defined as in the statement of Lemma
\ref{lem:Lemma 18} for $A$. Then there exists a direction $i$ for
which $\min\left(\left|A_{i,\,+}\right|,\,\left|A_{i,\,-}\right|\right)>\left|D_{+}\right|$. 
\end{lem}
\begin{proof}
Let $D_{+}$ and $D_{-}$ be defined as in the statement of Lemma
\ref{lem:Lemma 18}. For $I\subseteq X$, note that the set $A_{I}\subseteq Q_{n}$
defined by $A_{I}=\left\{ x\Delta I\,:\,x\in A\right\} $ is isomorphic
to $A$. By taking $I=\left\{ i\,:\,\left|A_{i,\,+}\right|>\left|A_{i,\,-}\right|\right\} $
it is easy to see that we may assume that $\left|A_{i,\,+}\right|\leq\left|A_{i,\,-}\right|$
holds for all $i$. Hence it suffices to show that under this condition
there exists $i$ for which $\left|A_{i,\,+}\right|>\left|D_{+}\right|$. 

Suppose that the claim is false. Then, by Lemma \ref{lem:Lemma 18},
we must have $\left|A_{i,\,+}\right|=\left|D_{+}\right|$ for all
$i$. Let $r$ be chosen such that $f_{r}<\left|A\right|<f_{r+1}$,
and note that such $r$ exists as $\left|A\right|\neq f_{s}$ for
all $s$. Since $N\left(B_{n-1}\left(\emptyset,r-1\right)\right)=B_{n-1}\left(\emptyset,r\right)$
and $f_{n-1,\,r-1}+f_{n-1,\,r}=f_{n,\,r}$, it follows from the maximality
assumption that $\left|D_{+}\right|\geq f_{n-1,\,r-1}$. On the other
hand, since $f_{n-1,\,r}+f_{n-1,\,r+1}=f_{n,\,r+1}$, the condition
$N\left(D_{+}\right)\subseteq D_{-}$ will be violated if $\left|D_{+}\right|\geq f_{n-1,\,r}$.
Combining these two bounds gives that 
\[
f_{n-1,\,r-1}\leq\left|D_{+}\right|<f_{n-1,\,r}.
\]
Hence there exists a set system ${\cal A}\subseteq X_{n}^{\left(r\right)}$
for which $D_{+}=X_{n}^{\left(\leq r-1\right)}\cup{\cal A}$ . Since
$N\left(D_{+}\right)\subseteq D_{-}$, it follows that $X_{n}^{\left(\leq r\right)}\cup\partial_{n}^{+}{\cal A}\subseteq D_{-}$. 

Let $C$ be the initial segment of the simplicial order of size $\left|A\right|$
on $Q_{n}$. Since $f_{r}<\left|A\right|<f_{r+1}$, it follows that
$C=X^{\left(\leq r\right)}\cup{\cal D}$ for some ${\cal D}\subseteq X^{\left(r+1\right)}$.
Combining $\left|A\right|=\left|D_{+}\right|+\left|D_{-}\right|$,
$D_{+}=X_{n}^{\left(\leq r-1\right)}\cup{\cal A}$ and $X_{n}^{\left(\leq r\right)}\cup\partial_{n}^{+}{\cal A}\subseteq D_{-}$,
it follows that $\left|{\cal D}\right|\geq\left|{\cal A}\right|+\left|\partial_{n}^{+}{\cal A}\right|$. 

For $B\subseteq Q_{n}$ define $f\left(B\right)=\sum_{x\in B}\left|x\right|$.
It is easy to see that among the sets $B\subseteq Q_{n}$ of given
size, $f\left(B\right)$ attains its minimum value when $B$ is taken
to be the initial segment of the simplicial order. Hence it follows
that $f\left(A\right)\geq f\left(C\right)$. It is easy to verify
that 
\begin{equation}
f\left(C\right)=\sum_{j=0}^{r}j{n \choose j}+\left(r+1\right)\left|{\cal D}\right|=n\sum_{j=0}^{r-1}{n-1 \choose j}+\left(r+1\right)\left|{\cal D}\right|.\label{eq:28}
\end{equation}
Let $\mathbb{I}\left\{ i\in x\right\} $ be the indicator function
of the event $\left\{ i\in x\right\} $ for $x\in Q_{n}$. Then for
any $B\subseteq Q_{n}$ we have 
\[
f\left(B\right)=\sum_{x\in B}\left|x\right|=\sum_{x\in B}\sum_{i=1}^{n}\mathbb{I}\left\{ i\in x\right\} =\sum_{i=1}^{n}\sum_{x\in B}\mathbb{I}\left\{ i\in x\right\} =\sum_{i=1}^{n}\left|B_{i,\,+}\right|.
\]
Since $\left|A_{i,\,+}\right|=\left|D_{+}\right|$ for all $i$, it
follows that 
\begin{equation}
f\left(A\right)=n\left|D_{+}\right|=n\sum_{j=0}^{r-1}{n-1 \choose j}+n\left|{\cal A}\right|.\label{eq:30}
\end{equation}
Since $f\left(A\right)\geq f\left(C\right)$, (\ref{eq:28}) and (\ref{eq:30})
imply that 
\begin{equation}
n\left|{\cal A}\right|\geq\left|{\cal D}\right|\left(r+1\right)\label{eq:31}
\end{equation}
Since $\left|{\cal D}\right|\ge\left|{\cal A}\right|+\left|\partial_{n}^{+}{\cal A}\right|$,
this implies that 
\begin{equation}
\left(n-r-1\right)\left|{\cal A}\right|\geq\left(r+1\right)\left|\partial_{n}^{+}{\cal A}\right|.\label{eq:32}
\end{equation}
The local LYM inequality \cite{key-18,key-1,key-2} for upper shadow
states that for ${\cal A}\subseteq X_{n}^{\left(r\right)}$ we have
\begin{equation}
\left(r+1\right)\left|\partial_{n}^{+}{\cal A}\right|\geq\left(\left(n-1\right)-r\right)\left|{\cal A}\right|,\label{eq:33}
\end{equation}
and the equality holds if and only if ${\cal A}=\emptyset$ or ${\cal A}=X_{n}^{\left(r\right)}$.
Combining (\ref{eq:32}) with (\ref{eq:33}), we obtain that the equality
must hold in (\ref{eq:33}). 

If ${\cal A}=\emptyset$, then (\ref{eq:31}) implies that ${\cal D}=\emptyset$.
Hence $C=X^{\left(\leq r\right)}$, which has size $f_{r}$. Hence
$\left|A\right|=f_{r}$, which contradicts the assumption. If ${\cal A}=X_{n}^{\left(r\right)}$,
then $D_{+}=X_{n}^{\left(\leq r\right)}$ and hence $X_{n}^{\left(\leq r+1\right)}=N\left(D_{+}\right)\subseteq D_{-}$.
Therefore $\left|A\right|=\left|D_{+}\right|+\left|D_{-}\right|\geq f_{r+1}$,
which contradicts the assumption that $\left|A\right|<f_{r+1}$. Thus
there must be a direction $i$ for which $\min\left(\left|A_{i,\,+}\right|,\,\left|A_{i,\,-}\right|\right)>\left|D_{+}\right|$.
\end{proof}
Now we are ready to prove Theorem \ref{thm:Theorem 4}. 

\addtocounter{thm}{-16}
\begin{thm}
Let $A\subseteq Q_{n}$
be a subset for which every $B\subseteq Q_{n}$ with $\left|B\right|=\left|A\right|$
satisfies $\left|N\left(B\right)\right|\geq\left|N\left(A\right)\right|$.
Then for all $t>0$ and $B\subseteq Q_{n}$ with $\left|B\right|=\left|A\right|$
we also have $\left|N^{t}\left(B\right)\right|\geq\left|N^{t}\left(A\right)\right|$.
\end{thm}
\addtocounter{thm}{15}
\begin{proof}
The proof is by induction on $n$. When $n\leq2$, it is easy to verify
that the result holds. Hence assume that $n\geq3$ and that the result
holds for $n-1$. 

Let $A\subseteq Q_{n}$ be a subset such that for all $B\subseteq Q_{n}$
with $\left|B\right|=\left|A\right|$ we have $\left|N\left(B\right)\right|\geq\left|N\left(A\right)\right|$, i.e.\ for which $\left|N(A)\right|$ is minimal.
By the argument presented in the proof of Lemma \ref{lem:Lemma 18},
we may assume that $\left|A_{i,\,+}\right|\leq\left|A_{i,\,-}\right|$
holds for all directions $i$. Let $D_{+}$ and $D_{-}$ be defined
as in the statement of Lemma \ref{lem:Lemma 18}. Then Lemma \ref{lem:Lemma 19}
implies that there exists a direction $i$ for which $\left|A_{i,\,+}\right|>\left|D_{+}\right|$.
For notational convenience, denote the $i$-sections in this direction
by $A_{+}$ and $A_{-}$. 

Let $C_{+}$ and $C_{-}$ be the initial segments of the simplicial
order on ${\cal P}\left(X_{i}\right)$ of the same sizes as $A_{+}$
and $A_{-}$ respectively, and set $C=C_{-}\cup\left(\left\{ i\right\} +C_{+}\right)$.
Then $\left|C_{+}\right|=\left|A_{+}\right|>\left|D_{+}\right|$,
and the maximality assumption on $D_{+}$ implies that $C_{-}\subseteq N\left(C_{+}\right)$.
Since $\left|C_{+}\right|\leq\left|C_{-}\right|$, it certainly follows
that $C_{+}\subseteq N\left(C_{-}\right)$. Since $N\left(A\right)$ is minimal and $C_{+}\subseteq N\left(C_{-}\right)$, 
Lemma \ref{lem:Lemma 10} with $t=1$ implies that $N\left(A_{+}\right)$ and $N\left(A_{-}\right)$
are minimal, and that

\begin{equation}
A_{\pm}\subseteq N\left(A_{\mp}\right).\label{eq:38}
\end{equation}
Since $N\left(A_{+}\right)$ and $N\left(A_{-}\right)$ are minimal,
the inductive hypothesis implies that for all sets $E_{+},\,E_{-}\subseteq{\cal P}\left(X_{i}\right)$
of the same sizes as $A_{+}$ and $A_{-}$ we have $\left|N^{t}\left(E_{+}\right)\right|\geq\left|N^{t}\left(A_{+}\right)\right|$
and $\left|N^{t}\left(E_{-}\right)\right|\geq\left|N^{t}\left(A_{-}\right)\right|$.
In particular, Harper's theorem implies that for all $t>0$ we have

\begin{equation}
\left|N^{t}\left(A_{\pm}\right)\right|=\left|N^{t}\left(C_{\pm}\right)\right|.\label{eq:39}
\end{equation}

By taking neighbourhood $t-1$ times from (\ref{eq:38}), it follows
that we have $N^{t-1}\left(A_{\pm}\right)\subseteq N^{t}\left(A_{\mp}\right)$
for all $t>0$. Since $C_{\pm}\subseteq N\left(C_{\pm}\right)$, we
also have $N^{t-1}\left(C_{\pm}\right)\subseteq N^{t}\left(C_{\mp}\right)$
for all $t>0$. Thus (\ref{eq:L2}) implies that 
\begin{equation}
\left|N^{t}\left(A\right)\right|=\left|N^{t}\left(A_{+}\right)\right|+\left|N^{t}\left(A_{-}\right)\right|\label{eq:34}
\end{equation}
and 
\begin{equation}
\left|N^{t}\left(C\right)\right|=\left|N^{t}\left(C_{+}\right)\right|+\left|N^{t}\left(C_{-}\right)\right|\label{eq:35}
\end{equation}
hold for all $t>0$. Combining (\ref{eq:39}), (\ref{eq:34}) and
(\ref{eq:35}), we obtain that $\left|N^{t}\left(A\right)\right|=\left|N^{t}\left(C\right)\right|$
for all $t>0$. Thus Harper's inequality implies that for all $B\subseteq Q_{n}$
with $\left|B\right|=\left|A\right|$ we have $\left|N^{t}\left(B\right)\right|\geq\left|N^{t}\left(A\right)\right|$,
which completes the proof of Theorem \ref{thm:Theorem 4}. 
\end{proof}
Now we can immediately apply Theorem \ref{thm:Theorem 4} to prove
that weak extremality implies extremality. 

\addtocounter{thm}{-17}
\begin{thm}
Let $A\subseteq Q_{n}$
be a weakly extremal subset. Then $A$ is extremal. 
\end{thm}
\addtocounter{thm}{16}
\begin{proof}
Since for all $B\subseteq Q_{n}$ with $\left|B\right|=\left|A\right|$
we have $\left|N\left(B\right)\right|\geq\left|N\left(A\right)\right|$,
Theorem \ref{thm:Theorem 4} implies that for all $B\subseteq Q_{n}$
with $\left|B\right|=\left|A\right|$ we have $\left|N^{t}\left(B\right)\right|\geq\left|N^{t}\left(A\right)\right|$.
Similarly applying Theorem \ref{thm:Theorem 4} to $A^{c}$ gives
that for all $B\subseteq Q_{n}$ with $\left|B\right|=\left|A\right|$
we have $\left|N^{t}\left(B^{c}\right)\right|\geq\left|N^{t}\left(A^{c}\right)\right|$.
Hence weak extremality implies extremality. 
\end{proof}
Theorem \ref{thm:Theorem 3} allows us to deduce that the classification
of weakly extremal sets can be done in the same way as the classification
of extremal sets. 
\begin{cor}
\label{cor:Corollary 20}Theorem \ref{thm:Theorem 2} holds with extremality
replaced by weak extremality. $\hfill\square$
\end{cor}

\section{A uniqueness result for certain sizes}

Recall that $f_{n,\,r}$ is defined to be the size of the exact Hamming
ball of radius $r$, and $g_{n,\,r}$ is defined to be the size of
the set $X^{\left(\leq r\right)}\cup\left(\left\{ 1\right\} +X_{1}^{\left(r\right)}\right)$.
In particular, recall that we have $f_{n,\,r}=\sum_{i=0}^{r}{n \choose i}$
and $g_{n,\,r}=f_{n,\,r}+{n-1 \choose r}$. 

Recall that in Section 2 we defined an extremal set $B_{r}$ of size
$g_{r}$ by setting $B_{r}=B\left(\emptyset,r\right)\cup B\left(\left\{ 1,2\right\} ,r\right)$.
The aim of this section is to prove that up to isomorphism the sets
$B_{r}$ introduced in Section 2 are the only sets of size $g_{r}$,
together with the initial segment, for which $N\left(A\right)$ is
minimal. Recently Keevash and Long \cite{key-16} studied stability
in the vertex isoperimetric inequality. Theorem \ref{thm:Theorem 21}
follows as a consequence of their more general result. 
\begin{thm}
\label{thm:Theorem 21}Let $r\leq n-1$, and let $A\subseteq Q_{n}$
with $\left|A\right|=g_{r}$ for which the size of $N\left(A\right)$
is minimal among the subsets of $Q_{n}$ of size $g_{r}$. Then either
$A$ is isomorphic to the initial segment of the simplicial order
or $A$ is isomorphic to $B_{r}$. 
\end{thm}
\begin{proof}
The main idea of the proof is to carefully analyse the codimension
1 compressions. Let $A\subseteq Q_{n}$ with $\left|A\right|=g_{r}$
for which $N\left(A\right)$ is minimal. As in the proof of Lemma
\ref{lem:Lemma 19}, by considering $A_{I}=\left\{ a\Delta I\,:\,a\in A\right\} $
if necessary for suitably chosen $I\subseteq X$, we may assume that
$\left|A_{i,\,+}\right|\leq\left|A_{i,\,-}\right|$ holds for all
directions $i$. 

Choose a direction $i$. Let $C_{i,\,+}$ and $C_{i,\,-}$ be the
initial segments of the simplicial order with $\left|C_{i,\,+}\right|=\left|A_{i,\,+}\right|$
and $\left|C_{i,\,-}\right|=\left|A_{i,\,-}\right|$, and define $C=C_{i,\,-}\cup\left(\left\{ i\right\} +C_{i,\,+}\right)$.
Let $D$ be the initial segment of the simplicial order of size $g_{r}$,
i.e.\ $D=X^{\left(\leq r\right)}\cup\left(\left\{ 1\right\} +X_{1}^{\left(r\right)}\right)$.
Recall that $N\left(D\right)=X^{\left(\leq r+1\right)}\cup\left(\left\{ 1\right\} +X_{1}^{\left(r+1\right)}\right)$
and that $\left|N\left(D\right)\right|=g_{r+1}$. Our first aim is
to find bounds for the sizes of $A_{i,\,-}$ and $A_{i,\,+}$. Since
their sizes are the same as the sizes of $C_{i,\,-}$ and $C_{i,\,+}$,
it suffices to prove these bounds for the sizes of $C_{i,\,-}$ and
$C_{i,\,+}$. 

If $\left|C_{i,\,-}\right|>g_{n-1,\,r}$, then we have $\left|C_{i,\,+}\right|<g_{n-1,\,r-1}$.
Hence we must have $\left|N\left(C_{i,\,+}\right)\right|\leq g_{n-1,\,r}<\left|C_{i,\,-}\right|$
and $\left|N\left(C_{i,\,-}\right)\right|\geq g_{n-1,\,r+1}$. Thus
it follows that $N\left(C_{i,\,+}\right)\subseteq C_{i,\,-}$, and
hence (\ref{eq:L2}) implies that 
\[
\left|N\left(C\right)\right|=\left|N\left(C_{i,\,-}\right)\right|+\left|C_{i,\,-}\right|>g_{n-1,\,r+1}+g_{n-1,\,r}=g_{n,\,r+1}.
\]
Since $\left|N\left(A_{i,\,-}\right)\right|\geq\left|N\left(C_{i,\,-}\right)\right|$
holds by Harper's inequality, (\ref{eq:L2}) implies that 
\begin{equation}
\left|N\left(A\right)\right|\geq\left|N\left(A_{i,\,-}\right)\right|+\left|A_{i,\,-}\right|\geq\left|N\left(C_{i,\,-}\right)\right|+\left|C_{i,\,-}\right|=\left|N\left(C\right)\right|>g_{n,\,r+1}.\label{eq:41}
\end{equation}
This contradicts the minimality of the size of $N\left(A\right)$.
Hence we must have $\left|C_{i,\,-}\right|\leq g_{n-1,\,r}$. 

Since $\left|C_{i,\,-}\right|\geq\left|C_{i,\,+}\right|$, it follows
that $\left|C_{i,\,-}\right|\geq\frac{1}{2}g_{n,\,r}=f_{n-1,\,r}$.
Combining these two gives that $f_{n-1,\,r}\leq\left|C_{i,\,-}\right|\leq g_{n-1,\,r}$.
Our next aim is to show that $\left|C_{i,\,-}\right|$, which equals
$\left|A_{i,\,-}\right|$, must always be either $f_{n-1,\,r}$ or
$g_{n-1,\,r}$. 
\addtocounter{thm}{-21}
\begin{claim}
For all $i$, we have $\left|C_{i,\,-}\right|=f_{n-1,\,r}$
or $\left|C_{i,\,-}\right|=g_{n-1,\,r}$.
\end{claim}
\addtocounter{thm}{20}
\begin{proof}
Since $f_{n-1,\,r}\leq\left|C_{i,\,-}\right|\leq g_{n-1,\,r}$ and
$\left|C\right|=g_{n,\,r}$, it follows that $g_{n-1,\,r-1}\leq\left|C_{i,\,+}\right|\leq f_{n-1,\,r}$.
In particular, it follows that we always have $C_{i,\,-}\subseteq N\left(C_{i,\,+}\right)$,
and we also have $C_{i,\,+}\subseteq N\left(C_{i,\,-}\right)$.

Let $j$ be the least element of $X_{i}$, i.e.\ $j=1$ if $i\neq1$
and $j=2$ if $i=1$. Since $C_{i,\,-}$ is an initial segment of
the simplicial order on ${\cal P}\left(X_{i}\right)$ with $f_{n-1,\,r}\leq\left|C_{i,\,-}\right|\leq g_{n-1,\,r}$
, it follows that $C_{i,\,-}=X_{i}^{\left(\leq r\right)}\cup\left(\left\{ j\right\} +{\cal A}\right)$
where ${\cal A}$ is an initial segment of the lexicographic order
on $X_{i,j}^{\left(r\right)}$. Similarly one can deduce that $C_{i,\,+}=X_{i}^{\left(\leq r-1\right)}\cup\left(\left\{ j\right\} +X_{i,j}^{\left(r-1\right)}\right)\cup{\cal B}$,
where ${\cal B}$ is an initial segment of the lexicographic order
on $X_{i,j}^{\left(r\right)}$. Note that 
\begin{align*}
 & \left|{\cal A}\right|+\left|{\cal B}\right|=\left|A\right|-f_{n-1,\,r-1}-f_{n-1,\,r}-\left|X_{i,j}^{\left(r-1\right)}\right|\nonumber \\
 & =g_{n,\,r}-f_{n,\,r}-{n-2 \choose r-1}={n-1 \choose r}-{n-2 \choose r-1}={n-2 \choose r}.
\end{align*}

It is easy to verify that 
\[
N\left(C_{i,\,-}\right)=X_{i}^{\left(\leq r+1\right)}\cup\left(\left\{ j\right\} +\partial_{i,j}^{+}{\cal A}\right)
\]
and 
\[
N\left(C_{i,\,+}\right)=X_{i}^{\left(\leq r\right)}\cup\left(\left\{ j\right\} +X_{i,j}^{\left(r\right)}\right)\cup\partial_{i,j}^{+}{\cal B},
\]
where again $\partial_{i,j}^{+}$ denotes the upper shadow operator
with respect to the ground set $X_{i,j}$. In particular, it follows
that 
\begin{equation}
\left|N\left(C_{i,\,-}\right)\right|=f_{n-1,\,r+1}+\left|\partial_{i,j}^{+}{\cal A}\right|\label{eq:43}
\end{equation}
and 
\begin{equation}
\left|N\left(C_{i,\,+}\right)\right|=g_{n-1,\,r}+\left|\partial_{i,j}^{+}{\cal B}\right|.\label{eq:44}
\end{equation}

Note that $f_{n-1,\,r+1}+g_{n-1,\,r}=g_{n,\,r+1}-{n-2 \choose r+1}$.
Since $C_{i,\,+}\subseteq N\left(C_{i,\,-}\right)$ and $C_{i,\,-}\subseteq N\left(C_{i,\,+}\right)$,
combining (\ref{eq:L2}) together with (\ref{eq:43}) and (\ref{eq:44})
we obtain that 
\begin{equation}
\left|N\left(C\right)\right|=\left|N\left(C_{i,\,+}\right)\right|+\left|N\left(C_{i,\,-}\right)\right|=g_{n,\,r+1}+\left|\partial_{i,j}^{+}{\cal A}\right|+\left|\partial_{i,j}^{+}{\cal B}\right|-{n-2 \choose r+1}.\label{eq:45}
\end{equation}
Applying The Local LYM inequality to ${\cal A}$ as a subset of $X_{i,j}^{\left(r\right)}$,
it follows that
\[
\left|\partial_{i,j}^{+}{\cal A}\right|\ge\frac{n-r-2}{r+1}\left|{\cal A}\right|,
\]
and the equality holds if and only if ${\cal A}_{i,j}=X^{(r)}$ or
${\cal A}=\emptyset$. Certainly the same holds for ${\cal B}$ as
well. Adding these together we obtain that 
\[
\left|\partial_{i,j}^{+}{\cal A}\right|+\left|\partial_{i,j}^{+}{\cal B}\right|\geq\frac{n-r-2}{r+1}\left(\left|{\cal A}\right|+\left|{\cal B}\right|\right).
\]
Since $\left|{\cal A}\right|+\left|{\cal B}\right|={n-2 \choose r}$
and $\frac{n-r-2}{r+1}{n-2 \choose r}={n-2 \choose r+1}$, it follows
that 
\begin{equation}
\left|\partial_{i,j}^{+}{\cal A}\right|+\left|\partial_{i,j}^{+}{\cal B}\right|\geq{n-2 \choose r+1}.\label{eq:47,5}
\end{equation}

Combining (\ref{eq:47,5}) with (\ref{eq:45}), we obtain that 
\begin{equation}
\left|N\left(C\right)\right|\geq g_{n,\,r+1}.\label{eq:48}
\end{equation}
Recall that $\left|N\left(A\right)\right|=g_{n,\,r+1}$. As in (\ref{eq:41}),
it follows from Harper's inequality that we have $\left|N\left(A\right)\right|\geq\left|N\left(C\right)\right|$,
and hence (\ref{eq:48}) implies that the equality holds in both applications
of the Local LYM inequality. In particular, we must have ${\cal A}=\emptyset$
or ${\cal A}=X{}_{i,j}^{\left(r\right)}$. In the first case we have
$\left|C_{i,\,-}\right|=f_{n-1,\,r}$ and in the second case we have
$\left|C_{i,\,-}\right|=g_{n-1,\,r}$. This completes the proof of
Claim 1. 
\end{proof}

If $\left|A_{i,\,-}\right|=f_{n-1,\,r}$ for some $i$, then we could
use Proposition \ref{prop:Prop5} to deduce that $A_{i,\,-}$ and
$A_{i,\,+}$ must be exact Hamming balls. The aim of the next claim
is to prove that such a direction $i$ must always exist. 

\addtocounter{thm}{-20}
\begin{claim}
There exists $i$ for which $\left|A_{i,\,-}\right|=f_{n-1,\,r}$.
\end{claim}
\addtocounter{thm}{19}

\begin{proof}
Suppose that the claim is false. Then by Claim 1 it
follows that $\left|A_{i,\,-}\right|=g_{n-1,\,r}$ for all $i$. As
in the proof of Lemma \ref{lem:Lemma 19}, for $B\subseteq Q_{n}$
define $f\left(B\right)=\sum_{x\in B}\left|x\right|$. Recall that
among the sets $B\subseteq Q_{n}$ of given size, $f\left(B\right)$
attains its minimum value when $B$ is taken to be the initial segment
of the simplicial order. Also recall that $f\left(B\right)=\sum_{i=1}^{n}\left|B_{i,\,+}\right|$
for any $B\subseteq Q_{n}$. 

Since $\left|A_{i,\,-}\right|=g_{n-1,\,r}$ for all $i$, it follows
that $\left|A_{i,\,+}\right|=g_{n-1,\,r-1}$ for all $i$. Hence $f\left(A\right)=ng_{n-1,\,r-1}$.
Let $D=X^{\left(\leq r\right)}\cup\left(\left\{ 1\right\} +X_{1}^{\left(r\right)}\right)$
be the initial segment of the simplicial order of size $g_{n,\,r}$.
It is easy to verify that $\left|D_{1,\,+}\right|=f_{n-1,\,r}$ and
$\left|D_{i,\,+}\right|=g_{n-1,\,r-1}$ for all $i\ge2$. Hence it
follows that $f\left(D\right)=f_{n-1,\,r}+\left(n-1\right)g_{n-1,\,r-1}$.
Since $D$ is the initial segment of the simplicial order of the same
size as $A$, it follows that $f\left(A\right)\geq f\left(D\right)$.
Thus we must have $g_{n-1,\,r-1}\geq f_{n-1,\,r}$, which is only
true if $r=n-1$. Since $g_{n-1,\,n-1}=2^{n}$, it follows that we
must have $A=Q_{n}$, in which case we also have  $\left|A_{i,\,-}\right|=2^{n-1}=f_{n-1,\,n-1}$,
as required. This completes the proof of Claim 2. 
\end{proof}

Let $i$ be a direction for which $\left|A_{i,\,-}\right|=f_{n-1,\,r}$,
and note that hence we also have $\left|A_{i,\,+}\right|=f_{n-1,\,r}$.
Hence it follows that $C_{i,\,+}\subseteq N\left(C_{i,\,-}\right)$
and $C_{i,\,-}\subseteq N\left(C_{i,\,+}\right)$, and thus Lemma
\ref{lem:Lemma 10} with $t=1$ implies that $N\left(A_{i,\,+}\right)$
and $N\left(A_{i,\,-}\right)$ are minimal. Since $\left|A_{i,\,+}\right|=\left|A_{i,\,-}\right|=f_{n-1,\,r}$,
Proposition \ref{prop:Prop5} implies that $A_{i,\,-}$ and $A_{i,\,+}$
are exact Hamming balls on ${\cal P}\left(X_{i}\right)$ with radius
$r$. Hence $A_{i,\,+}=B_{i}\left(x,r\right)$ and $A_{i,\,-}=B_{i}\left(y,r\right)$
for some $x,\,y\in Q_{n}$, and by symmetry we may assume that $x=\emptyset$. 

Note that (\ref{eq:L2}) implies that 
\begin{align*}
 & \left|N\left(A\right)\right|=\left|A_{i,\,+}\cup N\left(A_{i,\,-}\right)\right|+\left|A_{i,\,-}\cup N\left(A_{i,\,+}\right)\right|\\
 & =\left|B_{i}\left(\emptyset,r\right)\cup B_{i}\left(y,r+1\right)\right|+\left|B_{i}\left(\emptyset,r+1\right)\cup B_{i}\left(y,r\right)\right|.
\end{align*}
Recall that $\left|N\left(A\right)\right|=g_{n,\,r+1}=2f_{n-1,\,r+1}$
by the minimality of $\left|N\left(A\right)\right|$. Since $\left|B_{i}\left(y,r+1\right)\right|=\left|B_{i}\left(\emptyset,r+1\right)\right|=f_{n-1,\,r+1}$,
it follows that we must have $B_{i}\left(\emptyset,r\right)\subseteq B_{i}\left(y,r+1\right)$
and $B_{i}\left(y,r\right)\subseteq B_{i}\left(\emptyset,r+1\right)$.
In particular, we must have $d\left(y,\emptyset\right)\le1$. 

If $y=\emptyset$, it follows that $A$ is isomorphic to the initial
segment of the simplicial order under a map $\phi_{\sigma}$ for any
$\sigma\in S_{n}$ with $\sigma\left(i\right)=1$. If $y=\left\{ j\right\} $
for some $j\in X_{i}$, then $A$ is isomorphic to the set $B_{r}$
under a map $\phi_{\sigma}$ for any $\sigma\in S_{n}$ with $\sigma\left(i\right)=1$
and $\sigma\left(j\right)=2$. This completes the proof of Theorem
\ref{thm:Theorem 21}. 
\end{proof}


\begin{thebibliography}{10}
\bibitem{key-9}G. Aubrun, S.J. Szarek. Alice and Bob Meet Banach
- The Interface of Asymptotic Geometric Analysis and Quantum Information
Theory. Mathematical Surveys and Monographs, American Math. Soc, 2017.

\bibitem{key-10}B. Bollob\'{a}s. Combinatorics. Cambridge University
Press, 1986.

\bibitem{key-11}B. Bollob\'{a}s and I. Leader. Compressions and isoperimetric
inequalities. \textit{J. Comb. Theory Ser. A}, \textbf{56} (1): 47-62,
 1991.

\bibitem{key-12}Z. F\"{u}redi, J. R. Griggs. Families of finite sets
with minimum shadows. \textit{Combinatorica} \textbf{6} (4): 355-363,
1986.

\bibitem{key-13}L. Harper. Optimal numberings and isoperimetric problems
on graphs. \textit{J. Comb. Theory} \textbf{1} (3): 385-393, 1966.

\bibitem{key-14}G. Katona. The Hamming-sphere has minimum boundary.
\textit{Studia Sci. Math. Hungar. }\textbf{10}: 131-140, 1975.

\bibitem{key-15}G. Katona. A theorem of finite sets. \textit{Theory
of Graphs}, Akademiai Kiado, 187\textendash 207, 1968.

\bibitem{key-16}P. Keevash, E. Long. Stability for vertex isoperimetry
in the cube. arXiv preprint, arXiv:1807.09618. 

\bibitem{key-17}J. Kruskal. The number of simplices in a complex.
\textit{Mathematical Optimization Techniques}, University of California
Press, 251\textendash 278, 1963.

\bibitem{key-18}D. Lubell. A short proof of Sperner\textquoteright s
lemma. \textit{J. Comb. Theory}, \textbf{1} (2): 299, 1966. 

\bibitem{key-1}L. Meshalkin. Generalization of Sperner\textquoteright s
theorem on the number of subsets of a finite set. \textit{Theory of
Probability and Its Applications}, \textbf{8} (2): 203-204, 1963.

\bibitem{key-19}E. R\"{a}ty. Coordinate deletion of zeroes. \textit{Electron.
J. Combin. }\textbf{26} (3), Paper 3.50, 2019. 

\bibitem{key-2}K. Yamamoto. Logarithmic order of free distributive
lattice. \textit{Journal of the Mathematical Society of Japan}, \textbf{6}
(3-4): 343-353, 1954.
\end{thebibliography}
\end{document}